\magnification=\magstep1 
\def\mapright#1{\smash{\mathop{\longrightarrow}\limits^{#1}}} 

\newcount\refno\refno=0
\chardef\other=12
\newwrite\reffile
\immediate\openout\reffile=TempReferences
\outer\def\ref{\par\medbreak\global\advance\refno by 1
  \immediate\write\reffile{}
  \immediate\write\reffile{\noexpand\item{[\the\refno]}}
  \copytoblankline}
\def\copytoblankline{\begingroup\setupcopy\copyref}
\def\setupcopy{\def\do##1{\catcode`##1=\other}\dospecials
  \catcode`\|=\other \obeylines}
{\obeylines \gdef\copyref#1
  {\def\next{#1}%
  \ifx\next\empty\let\next=\endgroup %
  \else\immediate\write\reffile{\next} \let\next=\copyref\fi\next}}

\ref ALLEN, S.D.; SINCLAIR, A.M.; SMITH, R.R., The ideal structure of the Haagerup tensor product of C$^*$-algebras, {\sl J. reine angew. Math.}, 442 (1993) 111-148.

\newcount\ASS\ASS=\refno



\ref ARCHBOLD, R.J.; KANIUTH, E.; SCHLICHTING, G.; SOMERSET, D.W.B., Ideal spaces of the Haagerup tensor product of C$^*$-algebras, {\sl Internat. J. Math.}, 8 (1997) 1-29.

\newcount\AKSS\AKSS=\refno



\ref BECKHOFF, F., Topologies on the space of ideals of a Banach algebra, {\sl Stud. Math.}, 115 (1995) 189-205.

\newcount\Be\Be=\refno



\ref BECKHOFF, F.,  Topologies on the ideal space of a Banach algebra and spectral synthesis, {\sl Proc. Amer. Math. Soc.}, 125 (1997) 2859-2866.

\newcount\Beck\Beck=\refno

\ref BLECHER, D.P.,  Geometry of the tensor product of C$^*$-algebras, {\sl Math. Proc. Camb. Phil. Soc.}, 104 (1988) 119-127.

\newcount\Ble\Ble=\refno

\ref BLECHER, D.P., Tensor products which do not preserve operator algebras, {\sl Math. Proc. Camb. Phil. Soc.}, 108 (1990) 395-403.

\newcount\Blec\Blec=\refno

\ref BOHNENBLUST, H.F.; KARLIN, S., Geometrical properties of the unit sphere of Banach algebras, {\sl Ann. Math.}, (2) 62 (1955) 217-229.

\newcount\BK\BK=\refno

\ref BONSALL, F.F.; DUNCAN, J., {\sl Complete Normed Algebras}, Springer-Verlag, New York, 1973.

\newcount\BD\BD=\refno



\ref COHN, P.M., {\sl Algebra, Vol. 3} (2nd edn.), Wiley, Chichester, 1991.

\newcount\Cohn\Cohn=\refno





\ref FEINSTEIN, J.F.; SOMERSET, D.W.B., A note on ideal spaces of Banach algebras, {\sl Bull. London Math. Soc.}, 30 (1998) 611-617.

\newcount\FS\FS=\refno

\ref J.F. FEINSTEIN, D.W.B. SOMERSET, Strong regularity for uniform algebras, 
in {\sl Contemporary Math., 232}, `Proceedings of the 3rd Function Spaces Conference, Edwardsville, Illinois, May 1998', ed. K. Jarosz, pp. 139-149, Amer. Math. Soc., Rhode Island, 1999.

\newcount\Fein\Fein=\refno



\ref GIERZ, G.; HOFMANN, K.H.; KEIMEL, K.; LAWSON, J.; MISLOVE, M.;
SCOTT, D.S., {\sl A Compendium of Continuous Lattices}, Springer-Verlag,
New York, 1980.

\newcount\Comp\Comp=\refno

\ref GRAHAM, C.C.; McGEHEE, O.C.,\hfil{\sl Essays in Commutative 
Harmonic Analysis}, 
Springer-Verlag, New York, 1979.

\newcount\GM\GM=\refno



\ref HENRIKSEN, M.; KOPPERMAN, R.; MACK, J.; SOMERSET, D.W.B., Joincompact
spaces, continuous lattices, and C$^*$-algebras, {\sl Algebra Universalis}, 38 (1997) 289-323.

\newcount\HKMS\HKMS=\refno







\ref KAPLANSKY, I., The structure of certain operator algebras, {\sl Trans. Amer. Math. Soc.}, 70 (1951) 219-255.

\newcount\Kap\Kap=\refno

\ref KELLEY, J.L., {\sl General Topology}, Van Nostrand, Princeton, 1955.

\newcount\Kel\Kel=\refno

\ref KOPPERMAN, R., Asymmetry and duality in topology, {\sl Gen. Top. and Appl.}, 66 (1995) 1-39.

\newcount\Ko\Ko=\refno







\ref PALMER, T.W., {\sl Banach Algebras and the General Theory of $^*$-Algebras}, Vol. 1, C.U.P., New York, 1994.

\newcount\Pal\Pal=\refno

\ref PAULSEN, V.I.; SMITH, R.R., Multilinear maps and tensor norms on operator systems, {\sl J. Funct. Anal.} 73 (1987) 258-276.

\newcount\PS\PS=\refno

\ref RICKART, C.E., On spectral permanence for certain Banach algebras, {\sl Proc. Amer. Math. Soc.}, 4 (1953) 191-196.

\newcount\Rick\Rick=\refno

\ref RICKART, C.E., {\sl General Theory of Banach Algebras}, Van Nostrand, London, 1960.

\newcount\Ri\Ri=\refno

\ref RUDIN, W, {\sl Fourier Analysis on Groups}, Interscience, New York, 1962.

\newcount\Rud\Rud=\refno



\ref SOMERSET, D.W.B., Spectral synthesis for Banach algebras, {\sl Quart. J. Math. Oxford}, (2) 49 (1998) 501-521.

\newcount\Syn\Syn=\refno

\ref SOMERSET, D.W.B., Ideal spaces of Banach algebras, {\sl Proc. London Math. Soc.}, (3) 78 (1999) 369-400.

\newcount\Id\Id=\refno



\ref ZELAZKO, W., {\sl Banach Algebras}, Elsevier, Warsaw, 1973.

\newcount\Zel\Zel=\refno

\immediate\closeout\reffile

\centerline{\bf SPECTRAL SYNTHESIS FOR BANACH ALGEBRAS, II} 
\bigskip 
\centerline{\bf J. F. Feinstein and D. W. B. Somerset} 
\bigskip 
\centerline{Joel.Feinstein@nottingham.ac.uk}
\medskip
\centerline{ds@maths.abdn.ac.uk}
\bigskip 
\noindent{\bf Abstract} This paper continues the study of spectral synthesis 
and the topologies $\tau_{\infty}$ and $\tau_r$ on the ideal space of a Banach 
algebra, concentrating particularly on the class of Haagerup tensor products 
of C$^*$-algebras. For this class, it is shown that spectral synthesis is 
equivalent to the Hausdorffness of $\tau_{\infty}$. Under a weak extra 
condition, spectral synthesis is shown to be equivalent to the Hausdorffness 
of $\tau_r$. 
\bigskip 
\noindent{{\bf 1991 Maths Subject Classification} 46H10, 46K50} 
\bigskip 
\bigskip 
\noindent {\bf Introduction} 
\bigskip 
\noindent The notion of spectral synthesis is well-established for commutative 
Banach algebras and for $L^1$-group algebras. In [\the\Syn] the second author 
introduced a definition of spectral synthesis for a general unital Banach 
algebra. The motivation behind this was the discovery that a unital 
commutative Banach algebra $A$ has spectral synthesis if and only if the 
topology $\tau_{\infty}$, introduced by Beckhoff [\the\Be] on the set of 
closed ideals of $A$, is Hausdorff [\the\Syn; 2.6]. The definition of spectral 
synthesis introduced in [\the\Syn] was modelled on the properties of 
C$^*$-algebras, because $\tau_{\infty}$ was also known to be Hausdorff for 
this class [\the\Be]. The hope was that spectral synthesis might be equivalent 
to the Hausdorffness of $\tau_{\infty}$ for non-commutative Banach algebras. 
It was shown in [\the\Syn] that this was so for separable, unital PI-Banach 
algebras, and that in general the Hausdorffness of $\tau_{\infty}$ implies a 
weak form of spectral synthesis. Conversely a strong form of spectral 
synthesis implies that the topology $\tau_{\infty}$ is $T_1$. 

Because $\tau_{\infty}$ is seldom Hausdorff, the second author introduced 
another topology $\tau_r$ on the set $Id(A)$ of closed two-sided ideals of a 
Banach algebra $A$ [\the\Id]. This topology is always compact, like 
$\tau_{\infty}$, and it is Hausdorff whenever $\tau_{\infty}$ is Hausdorff 
[\the\Id; 3.1.1], and often even when $\tau_{\infty}$ is not. This is the 
case, for instance, with TAF-algebras [\the\Id], and with the algebra 
$C^1[0,1]$ [\the\FS]. On the other hand it was shown in [\the\FS] that for 
uniform algebras, $\tau_r$ is Hausdorff if and only if $\tau_{\infty}$ is 
Hausdorff. In [\the\Id] it was shown that if there is a compact Hausdorff 
topology on a subspace of $Id(A)$, which is related to the quotient norms in a 
useful way, then that topology necessarily coincides with the restriction of 
$\tau_r$. Thus for uniform algebras without spectral synthesis, such as the 
disc algebra, there is no useful compact Hausdorff topology on the space of 
closed ideals. 

This paper continues the study of the relationship between spectral synthesis 
and the Hausdorffness of $\tau_{\infty}$ and $\tau_r$. A slightly modified 
definition of spectral synthesis is introduced, for various reasons,
and several of the results of 
[\the\Syn] are extended to the non-unital situation. In particular it is shown 
that the problem of proving that the Hausdorffness of $\tau_{\infty}$ implies 
spectral synthesis is harder than one might imagine---one would first have to 
prove that there are no non-trivial algebraically simple Banach algebras. 

In the second part of the paper, we turn our attention to the class of Banach 
algebras obtained by taking the Haagerup tensor product of two C$^*$-algebras. 
There are various reasons for looking at this class. One is that spectral 
synthesis has already been studied for these algebras, and it is easy to find 
examples where synthesis holds, and others where it fails. A second reason is 
that the presence of the C$^*$-algebras makes the class reasonably tractable, 
and a considerable amount is already known about the ideal structure 
[\the\ASS], [\the\AKSS]. A third reason is that the Banach algebras in this 
class are in general neither Banach $^*$-algebras, nor operator algebras, so 
it might be hoped that the class is reasonably typical of semi-simple Banach 
algebras as a whole. We are able to show that spectral synthesis is equivalent 
to the Hausdorffness of $\tau_{\infty}$ for this class. 

In the final part of the paper we work with the same class of Banach algebras, 
this time considering the topology $\tau_r$. We show that if the 
C$^*$-algebras are unital and have the property that every closed prime ideal 
is maximal then their Haagerup tensor product has spectral synthesis if and 
only if the topology $\tau_r$ is Hausdorff. One novel feature of our approach 
is the use of the theory of continuous lattices. 

We now give the definitions of the various topologies on $Id(A)$, starting 
with the lower topology $\tau_w$. Let $A$ be a Banach algebra. A subbase for 
$\tau_w$ on $Id(A)$ is given by the sets $\{ I\in Id(A): I \not\supseteq J\}$ 
as $J$ varies through the elements of $Id(A)$. Thus the restriction of 
$\tau_w$ to the set of closed prime ideals is simply the hull-kernel topology. 
Next we define $\tau_{\infty}$. For each $k\in {\bf N}$, let $S_k=S_k(A)$ 
denote the set of seminorms (`seminorm' means `algebra seminorm' in this 
paper) $\rho$ on $A$ satisfying $\rho(a)\le k\Vert a\Vert$ for all $a\in A$. 
Then $S_k$ is a compact, Hausdorff space [\the\Be]. We say that 
$\rho\ge\sigma$, for $\rho,\sigma\in S_k$, if $\rho (a)\le\sigma (a)$ for all 
$a\in A$. The point of this upside-down definition is that if $\rho\ge\sigma$ 
then $\ker\rho\supseteq\ker\sigma$. Clearly if $\rho,\sigma\in S_k$ the 
seminorm $\rho\wedge\sigma$ defined by $(\rho\wedge\sigma)(a)=\max\, 
\{\rho(a),\sigma(a)\}$ is the greatest seminorm less than both $\rho$ and 
$\sigma$ in the order structure. Thus $S_k$ is a lattice. The topology 
$\tau_{\infty}$ is defined on $Id(A)$ as follows [\the\Be]: for each $k$ let 
$\kappa_k:S_k\to Id(A)$ be the map $\kappa_k(\rho)=\ker\rho$, and let $\tau_k$ 
be the quotient topology of $\kappa_k$ on $Id(A)$. Then 
$\tau_{\infty}=\bigcap_k\tau_k$. Clearly each $\tau_k$ is compact, so 
$\tau_{\infty}$ is compact. It is a useful fact that for $I\in Id(A)$, 
$Id(A/I)$ is $\tau_{\infty}$-- $\tau_{\infty}$ homeomorphic to the subset $\{ 
J\in Id(A):J\supseteq I\}$ of $Id(A)$ [\the\Be; Prop. 5]. 

Next we define the topology $\tau_r$, which is the join of two weaker 
topologies. The first is easily defined: $\tau_u$ is the weakest topology on 
$Id(A)$ for which all the norm functions $I\mapsto\Vert a+I\Vert$ $(a\in A,\ 
I\in Id(A))$ are upper semi-continuous. The other topology $\tau_n$ can be 
described in various different ways, but none is particularly easy to work 
with. A net $( I_{\alpha})$ in $Id(A)$ is said to have the {\sl normality 
property} with respect to $I\in Id(A)$ if $a\notin I$ implies that 
$\lim\inf\Vert a+I_{\alpha}\Vert >0$. Let $\tau_n$ be the topology whose 
closed sets $N$ have the property that if $(I_{\alpha})$ is a net in $N$ with 
the normality property relative to $I\in Id(A)$ then $I\in N$. It follows that 
if $(I_{\alpha})$ is a net in $Id(A)$ having the normality property relative 
to 
$I\in Id(A)$ then $I_{\alpha}\to I$ $(\tau_n)$. Any topology for which 
convergent nets have the normality property with respect to each of their 
limits (such a topology is said to have the {\sl normality property}) is 
necessarily stronger than $\tau_n$, but $\tau_n$ itself need not have the 
normality property. Indeed the following is true. Let $\tau_r$ be the topology 
on $Id(A)$ generated by $\tau_u$ and $\tau_n$. Then $\tau_r$ is always compact 
[\the\Id; 2.3], and $\tau_r$ is Hausdorff if and only if $\tau_n$ has the 
normality property [\the\Id; 2.12]. 
{ The topology $\tau_n$ is always stronger than 
$\tau_w$, and so $\tau_w$ and $\tau_n$ coincide on a given subset
of $Id(A)$ if
$\tau_w$ has the normality property on this subset.}

{ Finally, we introduce an auxiliary topology, $\tau_a$. This is
the topology generated by $\tau_u$ and $\tau_1$. Clearly we have 
$\tau_u \subseteq \tau_a$ and $\tau_\infty \subseteq \tau_a$.
By [\the\Id; 2.6], $\tau_n \subseteq \tau_\infty$ and so we also have
$\tau_r \subseteq \tau_a$. We shall need this topology in the proof of
Theorem 1.6.}

\bigskip 
\noindent The following simple lemma is taken from [\the\FS; 0.1]. 
\bigskip 
\noindent {\bf Lemma 0.1} {\sl Let $A$ be a Banach algebra. Let $(I_{\alpha})$ 
be a net in $Id(A)$, either 
decreasing or increasing, and correspondingly either set $I=\bigcap 
I_{\alpha}$ 
or $I=\overline{\bigcup I_{\alpha}}$. Then $I_{\alpha}\to I$ $(\tau_r)$.} 
\bigskip 
\bigskip 
\noindent {\bf 1. The foundations revisited} 
\bigskip 
\noindent In this section we re-examine the foundations for the work on 
spectral synthesis. We introduce a slightly more general definition of 
spectral synthesis and consider some of its elementary consequences. This 
change of definition obliges us to survey the results of [\the\Syn] to see how 
they are affected. Since we wish subsequently to work with group algebras, we 
also take the opportunity to liberate the theory from the requirement that the 
Banach algebra should be unital. This is sometimes more tricky than one 
anticipates. We consider some of the relations between spectral synthesis, 
weak spectral synthesis, $\tau_{\infty}$, and $\tau_r$ for commutative Banach 
algebras, PI-Banach algebras, and algebraically simple Banach algebras, and we 
also look at the stability properties of spectral synthesis on passage to 
ideals, quotients, and extensions. 
\bigskip 
\noindent{\bf Old definition of spectral synthesis} In [\the\Syn], a unital 
Banach algebra $A$ was said to have spectral synthesis if it has the following 
four properties: (i) $Prim(A)$ is locally compact (i.e. every point has a 
neighbourhood base of compact sets), (ii) every closed subset of $Prim(A)$ is 
a Baire space (i.e. the intersection of countably many dense open sets is 
dense), (iii) $\tau_w$ has the normality property on $Prim(A)$, (iv) $Id(A)$ 
is isomorphic to the lattice of open subsets of $Prim(A)$, under the 
correspondence $I\leftrightarrow \{ P\in Prim(A):P\not\supseteq I\}$. 

Here $Prim(A)$ is the space of primitive ideals of $A$ (i.e. kernels of 
algebraically irreducible representations with the hull-kernel topology). If 
the word $Prim(A)$ is replaced throughout by $Prime(A)$ (the set of proper 
closed prime ideals of $A$ with the hull-kernel topology) then $A$ is said to 
have {\sl weak spectral synthesis}, and if the word $Prim(A)$ is replaced 
throughout by $Max(A)$ (the set of closed, maximal modular ideals of $A$ with 
the hull-kernel topology) then $A$ is said to have {\sl strong spectral 
synthesis}. 
\bigskip 
Unfortunately, this old definition of spectral synthesis is slightly too 
restrictive for the algebras that we wish to consider. In these algebras, and 
also in Banach $^*$-algebras, it is natural to consider an ideal as 
`primitive' if it is the kernel of a topologically irreducible 
$*$-representation on a Hilbert space, and it is not evident that such an 
ideal is primitive in the usual algebraic sense. We shall therefore need to 
relax the requirements for spectral synthesis. Note, however, that we do not 
change the definitions for weak and strong spectral synthesis, except for 
dropping the requirement of an identity element. 

There is also a second reason for wanting to modify the definition of spectral 
synthesis. Recall that a non-empty closed subset of a topological space is 
{\sl irreducible} if it is not the union of two proper closed subsets. The 
closure of a point is a typical example of an irreducible closed set, and a 
topological space is said to be {\sl sober} if every irreducible closed set is 
the closure of a point. For instance, every Hausdorff space is sober. An 
infinite set with the cofinite topology is not sober, because it is an 
irreducible closed subset of itself, but is not the closure of any of its 
points. 

Sobriety of a space is closely connected with the Baire property. For example, 
a locally compact, sober space is a Baire space [\the\Comp; p.84]. 
Indeed if $X$ is a second countable $T_0$-space such that the sobrification 
$X^s$ (see below) is locally compact, then $X$ is sober if and only if every 
closed subspace of $X$ is a Baire space [\the\Comp; V.5.27(ii)]. This, of 
course, is axiom (ii) above, which was required for some of the arguments in 
[\the\Syn]. Thus axiom (ii) would be redundant if $Prim(A)$ were always sober. 
What we do, therefore, to obtain slightly more generality, is to replace 
$Prim(A)$ by its sobrification, as follows. 

Every $T_0$-space $X$ has a unique {\sl sobrification} $X^s$, which may be 
defined in the following way. Let $X^s$ be the topological space whose points 
are the irreducible closed sets of $X$, and whose non-empty open sets have the 
form $\{ A\in X^s: A\cap U\ne\emptyset\}$, where $U$ varies through the 
non-empty open subsets of $X$. Then $X^s$ is a sober space. The map $i:X\to 
X^s$, $i(x)=\{ x\}^-$, embeds $X$ in $X^s$. Evidently $X$ and $X^s$ have 
isomorphic lattices of open sets, and $i(X)= X^s$ if and only if $X$ is sober. 
\bigskip 
\noindent For an algebra $R$ over ${\bf C}$, let $Idl(R)$ denote the set of 
all two-sided ideals of $R$. An ideal $I\in Idl(R)$ is said to be {\sl 
semisimple} if it is an intersection of primitive ideals, and to be {\sl 
strongly semisimple} if it is an intersection of maximal modular ideals. If 
$R$ is a Banach algebra then semisimple and strongly semisimple ideals are, of 
course, automatically closed. 

\bigskip 
\noindent{\bf Lemma 1.1} {\sl Let $R$ be an algebra over ${\bf C}$. The 
sobrification $Prim^s(R)$ of $Prim(R)$ is (homeomorphic to) the set of 
semisimple, prime ideals of $R$, with the hull-kernel topology. The 
sobrification $Max^s(R)$ of $Max(R)$ is (homeomorphic to) the set of strongly 
semisimple, prime ideals of $R$, with the hull-kernel topology.} 
\bigskip 
\noindent {\bf Proof.} A simple argument shows that for a semisimple ideal 
$P$, the hull-kernel closed set $X=\{ Q\in Prim(R):Q\supseteq P\}$ is 
irreducible if and only if $P$ is prime. Thus there is a one-to-one 
correspondence between $Prim^s(R)$ and the set of semisimple prime ideals, and 
it is straightforward to confirm that the map is a homeomorphism. 

An analogous argument deals with $Max(R)$. Q.E.D. 
\bigskip 
\noindent For example, if $A$ is the disc algebra then $Prim(A)=Max(A)$ is an 
irreducible closed subset of itself in the hull-kernel topology, but is not 
the closure of any of its points. In fact $Prim^s(A)=Max(A)\cup\{\,\{ 
0\}\,\}$. On the other hand, since every Hausdorff space is sober, 
$Prim^s(A)=Max(A)$ when $A$ is a completely regular, commutative Banach 
algebra. If $A$ is a C$^*$-algebra then every closed ideal of $A$ is 
semisimple, so $Prim^s(A)=Prime(A)$. The famous open problem of whether every 
closed, prime ideal of a (non-separable) C$^*$-algebra $A$ is primitive is 
thus precisely the question whether $Prim(A)$ is a sober space for an 
arbitrary C$^*$-algebra. 
\bigskip 
\noindent We are now ready to introduce the new definition of spectral 
synthesis. 
\bigskip 
\noindent {\bf Definition 1.2} A Banach algebra $A$, unital or otherwise, has 
{\sl spectral synthesis} if it has the following properties: 

(i)$'$, $Prim^s(A)$ (the space of semisimple, prime ideals with the 
hull-kernel topology) is locally compact, 

(ii)$'$ $\tau_w$ has the normality property on $Prim(A)$, 

(iii)$'$ $Id(A)$ is isomorphic to the lattice of open subsets of $Prim(A)$, 
under the correspondence $I\leftrightarrow \{ P\in Prim(A):P\not\supseteq 
I\}$. 
Equivalently, every proper, closed ideal of $A$ is semisimple. 

\bigskip 
\noindent A few remarks on this new definition are in order. 
\smallskip 
\noindent {\bf Remarks} (a) It can be shown [\the\Comp; V.5.10] that if $X$ is 
a $T_0$-space and $Y$ is locally compact with $X\subseteq Y\subseteq X^s$ then 
$X^s$ is locally compact. Thus if $Prim(A)$ is locally compact, 
{ then (since $Prim(A)$ is a $T_0$-space)}
$Prim^s(A)$ is 
also locally compact. Every Banach algebra, therefore, which had spectral 
synthesis under the old definition still has it under the new definition. As a 
matter of fact, the authors do not know of a unital Banach algebra which has 
spectral synthesis under the new definition but not under the old. 

(b) $Prim^s(A)$ is always sober, so if it is also locally compact then every 
closed subspace of $Prim^s(A)$ is a Baire space [\the\Comp; V.5.27(ii)]. Thus 
we do not need any axiom corresponding to (ii) in the old definition of 
spectral synthesis. 

(c) Axiom (ii)$'$ is, of course, axiom (iii) from the old definition of 
spectral synthesis. It is easy to show that $\tau_w$ has the normality 
property on $Prim(A)$ if and only if it has the normality property on 
$Prim^s(A)$. 

(d) It was shown in [\the\Syn; 1.1] that if $A$ has weak spectral synthesis 
then $\tau_w$ has the normality property on $Id(A)$. The proof works perfectly 
well in the non-unital case too. Hence it follows from Proposition 1.10 
(below) that if $A$ has spectral synthesis (or strong spectral synthesis) then 
$\tau_w$ has the normality property on $Id(A)$. This implies that $\tau_w$ and 
$\tau_n$ coincide on $Id(A)$ [\the\Id; p.373 and p.375]. 

(e) Axiom (iii)$'$ is, of course, axiom (iv) from the old definition of 
spectral synthesis. Recall that $Prim(A)$ and $Prim^s(A)$ have isomorphic 
lattices of open sets. 

\bigskip 
\noindent The next thing to do is to check that the new definition coincides 
with the standard one for commutative Banach algebras. Recall that a (possibly 
non-unital) commutative Banach algebra $A$ has spectral synthesis (usual 
definition) if the map $I\mapsto \{P\in Prim(A):P\supseteq I\}$ sets up a 1--1 
correspondence between closed ideals of $A$ and Gelfand closed subsets of 
$Prim(A)$. This is equivalent to requiring that the hull-kernel and Gelfand 
topologies coincide on $Prim(A)$, and that every closed ideal of $A$ is 
semisimple. 

Let $A$ be a commutative Banach algebra. Then every primitive ideal of $A$ is 
maximal and modular, indeed the kernel of a character [\the\Pal; 4.2.19], but 
if $A$ is non-unital then $A$ might have maximal proper ideals which are not 
primitive. Here is an instructive example from [\the\Beck; Example 2]. 
\bigskip 
\noindent {\bf Example 1.3} Let $A=C[0,1]$ be the commutative C$^*$-algebra of 
continuous functions on the interval $[0,1]$ with the supremum norm. Let $z\in 
A$ be the identity map given by $z(t)=t$ $(t\in[0,1])$. Define a new 
multiplication $\diamond$ on $A$ by $f\diamond g=fzg$, and let $B$ be the 
resulting Banach algebra. Then $Id(B)=Id(A)$, and the character space of $B$ 
is the set $\{ t\,\omega_t: t\in (0,1]\}$ where $\omega_t$ is the point 
evaluation at $t$. The ideal $\ker\omega_0$ is a maximal ideal of $B$, but not 
a primitive ideal, and $B/\ker\omega_0$ is isomorphic to the complex numbers 
with the zero multiplication. 

The next lemma (for which we have not been able to find a reference) shows 
that this example is typical. 
\bigskip 
\noindent {\bf Lemma 1.4} {\sl Let $A$ be a commutative Banach algebra and let 
$M$ be a maximal proper ideal of $A$. Then $A/M$ is isomorphic to ${\bf C}$, 
either with the usual multiplication or with the zero multiplication.} 
\bigskip 
\noindent {\bf Proof.} Set $B=A/M$. Let $b\in B\setminus\{0\}$. Then $Bb$ is 
an ideal of $B$, so either $Bb=\{0\}$ or $Bb=B$. In the first case, $\{\lambda 
b:\lambda\in{\bf C}\}$ is a non-zero ideal of $B$, so $\{\lambda 
b:\lambda\in{\bf C}\}=B$. Hence $B$ is one-dimensional with zero 
multiplication. In the second case, there exists $u\in B$ such that $ub=b$. 
Since $Bb=B$, for any $c\in B$ there exists $d\in B$ such that $db=c$, and 
then $uc=udb=db=c$. Thus $u$ is the identity for $B$. Hence $M$ is a maximal 
modular ideal of $A$, so $M$ is closed and $B=\{ \lambda u:\lambda \in {\bf 
C}\}$, see [\the\BD; \S 16 Theorem 5]. Q.E.D. 
\bigskip 
\noindent If $M$ above is not modular then $M$ need not be closed. For 
example, let $A$ be an infinite-dimensional Banach space with the zero 
multiplication, and let $M$ be a dense hyperplane. Then $A$ is a commutative 
Banach algebra, and $M$ is a maximal ideal of $A$ which is not closed. 

If $A$ is a commutative Banach algebra and $M$ is a maximal modular ideal of 
$A$ then, for $a\in A$, we shall identify the coset $a+M$ with the value that 
the character corresponding to $M$ takes at $a$. One inconvenience with 
non-unital commutative Banach algebras is that the process of evaluating at a 
character and taking the modulus sometimes gives a number strictly less than 
the quotient norm for the corresponding maximal modular ideal. This happens in 
Example 1.3, for instance. To get round this problem, we use the following 
lemma. 
\bigskip 
\noindent{\bf Lemma 1.5} {\sl Let $A$ be a commutative Banach algebra. Suppose 
that $(M_{\alpha})$ is a net in $Prim(A)$ and that $a\in A$ with 
$\lim_{\alpha} a+M_{\alpha}= 0$. Then $\lim_{\alpha} \Vert 
a^2+M_{\alpha}\Vert= 0$.} 
\bigskip 
\noindent{\bf Proof.} For any character $\phi$ on $A$, 
$a^2-\phi(a)a\in\ker\phi$. Hence $$\Vert a^2+\ker\phi\Vert\le\Vert 
a^2-(a^2-\phi(a)a)\Vert=\Vert \phi(a) a\Vert=\vert\phi(a)\vert\,\Vert 
a\Vert.$$ 
Thus $\limsup_{\alpha}\Vert a^2+M_{\alpha}\Vert\le \limsup_{\alpha}\vert 
a+M_{\alpha}\vert\,\Vert a\Vert =0$. Q.E.D. 
\bigskip 
\noindent For a Banach algebra $A$, let ${\cal M}'$ be the space of maximal 
modular ideals of codimension one. The next result was proved in [\the\Be; 
Proposition 11] under the additional requirement that $A$ has a bounded 
approximate identity (see also Maths Reviews 97f:46073, where it is observed 
that one only need assume that $A^2=A$). Part of [\the\Be; Proposition 11] was 
that the set ${\cal M}'\cup\{ A\}$ is $\tau_{\infty}$-closed if $A$ has a 
bounded approximate identity. In Example 1.3, however, $\ker \omega_0$ is in 
the $\tau_{\infty}$-closure of ${\cal M}'$, so that particular part of 
[\the\Be; Proposition 11] cannot be extended to the general case. 
\bigskip 
\noindent {\bf Theorem 1.6} {\sl Let $A$ be a Banach algebra. Then the 
restrictions of the topologies $\tau_n$, $\tau_r$, and $\tau_{\infty}$ all 
coincide on the Gelfand space ${\cal M}'$ of $A$ with the Gelfand topology.} 
\bigskip 
\noindent {\bf Proof.} The set ${\cal M}'$ is $\tau_w$-closed in $Prim(A)$ by 
a result of Kaplansky's, see [\the\Id; 5.1]. Set $K=\bigcap \{M:M\in {\cal 
M}'\}$. Then $A/K$ is commutative, and ${\cal M}'$ is homeomorphic, in all 
four topologies under consideration, to the Gelfand space of $A/K$ [\the\Be; 
Proposition 5], [\the\Id; 2.9]. Thus we may assume at the outset that $A$ is 
commutative. Under this assumption we have 
${\cal M}'=Prim(A)$, so to complete the proof
we need to show that the topologies under 
consideration agree on $Prim(A)$.

Let $X_1$ be the set of closed ideals of $A$ of codimension not greater than 
one. Then $X_1$ is a $\tau_n$-closed subset of $Id(A)$ [\the\Id; 5.1].
{ By 
[\the\Id; 4.3(ii)] the topologies $\tau_n$ and $\tau_a$ coincide on 
the set of minimal elements of $X_1$, which is $X_1\setminus\{ A\}$
(recall that $\tau_a$ is the topology generated by $\tau_u$ and
$\tau_1$). However we know that $\tau_n \subseteq \tau_\infty
\subseteq \tau_a$ and also $\tau_n \subseteq \tau_r \subseteq \tau_a$. 
It follows that all four of these topologies coincide on 
$X_1\setminus\{ A\}$, and hence that $\tau_n$, $\tau_r$ and 
$\tau_\infty$ coincide on the subset $Prim(A)\subseteq 
X_1\setminus\{ A\}$.}

Suppose that $(P_{\alpha})$ is a net in $Prim(A)$ converging to $P\in Prim(A)$ 
in the Gelfand topology. For each $\alpha$ let $c_{\alpha}$ be the seminorm in 
$S_1(A)$ given by $c_{\alpha}(a)=\vert a+P_{\alpha}\vert$ $(a\in A)$. Then 
$c_{\alpha}\to c$ where $c\in S_1(A)$ is the seminorm $c(a)=\vert a+P\vert$ 
$(a\in A)$. Hence $P_{\alpha}=\ker c_{\alpha}\to\ker c=P$ in the 
$\tau_{\infty}$ topology. Thus the Gelfand topology on $Prim(A)$ is stronger 
than the relative $\tau_{\infty}$ topology. 

Conversely, let $Y$ be a Gelfand-open subset of $Prim(A)$ with compact closure 
(in the Gelfand topology). We shall show that $Z=X_1\setminus Y$ is 
$\tau_n$-closed in $Id(A)$. This will show that $Y$, and hence every 
Gelfand-open subset of $Prim(A)$, is $\tau_n$-open in $Prim(A)$. It follows 
that the various topologies coincide on $Prim(A)$. 

Let $(I_{\alpha})$ be a net in $Z$ and let $I\in Id(A)\setminus Z$. We have to 
show, by [\the\Id; 2.5], that there exists $a\in A\setminus I$ such that 
$\lim\inf\Vert a+I_{\alpha}\Vert=0$. If $I\notin X_1$ then this follows at 
once from [\the\Id; 2.5] and the fact that $X_1$ itself is $\tau_n$-closed in 
$Id(A)$. So suppose that $I\in Y$. If $(I_{\alpha})$ is eventually in 
$X_1\setminus Prim(A)$ then choose any $a\in A\setminus I$. Then $a^2\notin 
I$, but $a^2\in I_{\alpha}$ eventually by Lemma 1.4. Hence $\lim_{\alpha}\Vert 
a^2+I_{\alpha}\Vert=0$. 
The other possibility is that $(I_{\alpha})$ is frequently in $Prim(A)$. The 
local compactness of the Gelfand topology on $Prim(A)$ implies that 
$(I_{\alpha})$ has a subnet $(I_{\beta})$ in $Prim(A)$ such that either 
$(I_{\beta})$ goes to infinity (i.e. is eventually outside every Gelfand 
compact subset of $Prim(A)$) or for which there exists $J\in Prim(A)$ 
with $I_{\beta}\to J$ in the Gelfand topology. In the first case choose any 
$a\notin I$; in the second case choose $a\in J\setminus I$. Then in either 
case $a\notin I$ but $a+I_{\beta}\to 0$. Hence $\lim_{\beta}\Vert 
a^2+I_{\beta}\Vert=0$ by Lemma 1.5, but $a^2\notin I$, since $I\in Prim(A)$. 
Q.E.D. 
\bigskip 
\noindent {\bf Corollary 1.7} {\sl Let $A$ be a commutative Banach algebra. 
Then $A$ has spectral synthesis (in the sense of this paper) if and only if 
$A$ has spectral synthesis in the usual sense.} 
\bigskip 
\noindent{\bf Proof.} Suppose that $A$ has spectral synthesis in the sense of 
this paper. Then, { since 
$\tau_w$ has the normality property on $Prim(A)$, $\tau_w$
and $\tau_n$ coincide on $Prim(A)$.}
Thus the hull-kernel { 
topology ($\tau_w$) and the Gelfand topology} coincide on $Prim(A)$ by Theorem 
1.6. Condition (iii)$'$, on the other hand, implies that every closed ideal of 
$A$ is semisimple. Hence $A$ has spectral synthesis in the usual sense. 

Conversely, suppose that $A$ has spectral synthesis in the usual sense. Then 
certainly condition (iii)$'$ holds. Furthermore, since $Prim(A)$ is 
$\tau_w$-Hausdorff it is sober, so $Prim(A)=Prim^s(A)$. Thus condition (i)$'$ 
holds. Finally condition (ii)$'$ follows from the fact that for any $a\in A$ 
and $P\in Prim(A)$, $\Vert a+P\Vert\ge\vert a+P\vert$, and the function 
$P\mapsto |a+P|$ is $\tau_w$-continuous on $Prim(A)$ since $A$ is completely 
regular. Thus $A$ has spectral synthesis in the sense of this paper. Q.E.D. 
\bigskip 
\noindent The next result was established for unital, commutative Banach 
algebras 
in [\the\Syn; 2.6]. We follow exactly the same method of proof, but making use 
of Theorem 1.6. Recall the theorem of Rickart's [\the\Rick] that if $A$ is a 
completely regular, semisimple, commutative Banach algebra then every norm 
$\Vert\,.\,\Vert'$ on $A$ is spectral, that is, for all $a\in A$, $\Vert 
a\Vert'\ge\sup\,\{\vert a+M\vert:M\in Max(A)\}$. This implies that $\{ 0\}$ is 
$\tau_{\infty}$-closed in $Id(A)$. Recall also the theorem of Bohnenblust and 
Karlin [\the\BK], see [\the\Zel; 12.7], that if $A$ is a commutative Banach 
algebra and $a\in A$ then the spectral radius of $a$ is the infimum of $\Vert 
a\Vert'$ over all possible norms $\Vert\, .\, \Vert'$ equivalent to the 
original norm, and bounded by it. Thus if $a$ is a non-zero quasinilpotent 
element, there is a sequence $(\Vert\, .\, \Vert_n)$ of norms in $S_1(A)$ such 
that $\Vert a\Vert_n\to 0$. { Considering a cluster point of such
a sequence we see that $\{0\}$ is not $\tau_{\infty}$-closed in $Id(A)$.}
It follows that a necessary condition for $\{ 0\}$ 
to be $\tau_{\infty}$-closed in $Id(A)$ is that $A$ should be semisimple. 
\bigskip 
\noindent{\bf Theorem 1.8} {\sl Let $A$ be a commutative Banach algebra. Then 
$A$ has spectral synthesis if and only if the topology $\tau_{\infty}$ is 
Hausdorff on $Id(A)$.} 
\bigskip 
\noindent {\bf Proof.} Suppose that $\tau_{\infty}$ is Hausdorff on $Id(A)$. 
Then $\tau_{\infty}$ is Hausdorff on $Id(A/I)$ for all $I\in Id(A)$ [\the\Be; 
Prop. 5], so $A/I$ must be semisimple, for all $I\in Id(A)$ by the theorem of 
Bohnenblust and Karlin just mentioned. Thus condition (iii)$'$ holds. 
Furthermore [\the\Syn; Theorem 2.5] shows that $\tau_w$ has the normality 
property on $Id(A)$, hence on $Prim(A)$, so condition (ii)$'$ holds. It also 
follows that $\tau_w$ and $\tau_n$ are equal, so Theorem 1.6 shows that the 
Gelfand and hull-kernel topologies coincide on $Prim(A)$. Hence condition 
(i)$'$ holds, and $A$ has spectral synthesis. 

{ The proof of the converse is identical to the corresponding
part of the proof of [\the\Syn; 2.6]; there are no difficulties in 
passing to the non-unital case.}
Q.E.D. 
\bigskip 
\noindent In general a commutative Banach algebra can have $\tau_r$ Hausdorff 
and yet fail to have spectral synthesis. Various examples are given in 
[\the\Id; Section 3]. For the class of uniform algebras, however, it was shown 
in [\the\FS; 1.2] that $\tau_r$ cannot be Hausdorff unless spectral synthesis 
holds. Our next theorem is a more general version of that result, valid also 
for non-unital Banach algebras. 
\bigskip 
\noindent {\bf Definition} Let $A$ be a Banach function algebra, and 
let $\Gamma(A)$ denote the Shilov boundary of $A$. Recall that a Gelfand 
compact subset $X$ of $Max(A)$ is a {\sl Helson set} if $A|_X=C(X)$ (where 
$C(X)$ is the algebra of continuous complex functions on $X$). Letting $I$ be 
the closed ideal consisting of elements of $A$ which vanish on $X$, the least 
constant $K$ such that $$K\sup\{|f(x)|:x\in X\}\ge \Vert f+I\Vert\ \hbox{ for 
all }f\in A$$ is called the {\sl Helson constant} of $X$. We say that a Banach 
function algebra $A$ has the {\sl Helson property} (with constant $k$) if 
there is a constant $k$ such that whenever $U$ is a non-empty Gelfand open 
subset of $\Gamma(A)$ there is an increasing net $(F_{\alpha})_{\alpha}$ of 
Helson sets of constant bounded by $k$ contained in $U$ such that 
$\bigcup_{\alpha} F_{\alpha}$ is Gelfand dense in $U$. 

For example, if $A$ is a uniform algebra then the collection of finite p-sets 
in the open set $U$ is an increasing net of Helson sets, with Helson constant 
$1$, whose union is dense in $U$, see [\the\FS]. Thus uniform algebras have 
the Helson property with constant $1$. We shall see after Proposition 3.5 that 
if $A$ and $B$ are commutative C$^*$-algebras then the Haagerup tensor product 
$A\otimes_h B$ (when viewed as a function algebra on its maximal ideal space) 
also has the Helson property with constant $1$. 

For a Banach function algebra $A$, and a Gelfand closed subset $F$ of 
$Max(A)$, let $I(F)$ be the ideal of elements of $A$ which vanish on $F$, and 
let $J(F)$ be the ideal of elements of compact support vanishing in a Gelfand 
neighbourhood of $F$ in $Max(A)$. If $F\subseteq\Gamma(A)$, let $L(F)$ be 
the ideal obtained as the closure of the set of elements of $A$ having compact 
support on $\Gamma(A)$ and vanishing in a Gelfand neighbourhood 
of $F$ in $\Gamma(A)$. If $Max(A)=\Gamma(A)$ then of course $L(F)$ is simply 
the 
closure of $J(F)$, but this need not be true when $Max(A)\ne \Gamma(A)$. 

Now recall that one characterization of spectral synthesis for a Banach 
function algebra $A$ is that $A$ has spectral synthesis if $J(F)$ is dense in 
$I(F)$ for each Gelfand closed subset $F$ of $Max(A)$. A related, weaker 
notion is that $A$ is {\sl strongly regular} if $J(\{ x\})$ is dense in $I(\{ 
x\})$ for each $x\in Max(A)$. If $I(\{ x\})=L(\{ x\})$ for each $x\in 
\Gamma(A)$ then $A$ is {\sl strongly regular on $\Gamma(A)$}; in fact, this 
implies that $\Gamma(A)=Max(A)$, see [\the\Fein; Theorem 2] (the simple 
argument there is for unital Banach function algebras, but it is easily 
modified to cope with the non-unital case). 
\bigskip 
\noindent {\bf Theorem 1.9} {\sl Let $A$ be a Banach function algebra with the 
Helson property (with constant $k$). Then $A$ has spectral synthesis if and 
only if $\tau_r$ is Hausdorff on $Id(A)$.} 
\bigskip 
\noindent {\bf Proof.} If $A$ has spectral synthesis then $\tau_r$ is 
Hausdorff by Theorem 1.8 and [\the\Id; 3.1.1]. If $A$ does not have spectral 
synthesis, there are two possibilities. Either 
$Max(A)=\Gamma(A)$, in which case $J(F)$ is dense in $L(F)$ for every closed 
subset $F$ of $Max(A)$, so by assumption there is a Gelfand closed subset $X$ 
of 
$Max(A)$ such that $I(X)\ne L(X)$. Otherwise $Max(A)\ne \Gamma (A)$, so 
$A$ is not strongly regular on $\Gamma (A)$, as we remarked just above. Thus 
there exists $x\in \Gamma (A)$ such that $I(\{ x\})\ne L(\{ x\})$. Hence in 
either case there is a Gelfand closed subset $X$ of $\Gamma(A)$ such that 
$I(X)\ne L(X)$. 

Let $(V_{\alpha})_{\alpha}$ be a net of decreasing, open neighbourhoods of 
$X$ in $\Gamma(A)$, each having compact complement in $\Gamma(A)$, such that 
$\bigcap_{\alpha} N_{\alpha}=X$ (where for each $\alpha$, $N_{\alpha}$ is the 
closure of $V_{\alpha}$ in the Gelfand topology). Then 
$(I(N_{\alpha}))_{\alpha}$ is an increasing net in $Id(A)$, and 
$I(N_{\alpha})\subseteq L(X)$, for each 
$\alpha$, so $$I:=\overline{\bigcup_{\alpha} I(N_{\alpha})}\subseteq L(X).$$ 
For each $\alpha$, let $(F_{\beta(\alpha)})_{\beta(\alpha)}$ be an 
increasing net of Helson sets in $V_{\alpha}$ of constant bounded by $k$
and whose union is dense in $V_{\alpha}$.
Then 
$(I({F_{\beta(\alpha)}}))_{\beta(\alpha)}$ is a decreasing net in 
$Id(A)$, and the density condition implies that 
$\bigcap_{\beta(\alpha)}I({F_{\beta(\alpha)}}) = 
I(N_{\alpha})$. Hence $I({F_{\beta(\alpha)}})\mapright{\beta(\alpha)} 
I(N_{\alpha})$ $(\tau_r)$ by Lemma 0.1. But $I(N_{\alpha})\to 
I$ $(\tau_r)$, also by Lemma 0.1, so if $(I(F_{\gamma}))_{\gamma}$ denotes the 
`diagonal' net, see [\the\Kel; \S 2, Theorem 4], then $I(F_{\gamma})\to 
I\subseteq L(X)$ $(\tau_r)$. 

Suppose that $f\notin I(X)$. Then there is a Gelfand open subset $U$ of 
$\Gamma(A)$ 
meeting $X$, and an $\epsilon >0$ such that $\vert f(x)\vert>\epsilon$ for all 
$x\in U$. By the density condition there is, for each 
$\alpha$, a $\beta_0(\alpha)$ such that 
$F_{\beta_0(\alpha)}\cap U$ is non-empty, and hence such that 
$\Vert f+I(F_{\beta(\alpha)})\Vert>\epsilon$ for all 
$\beta(\alpha)\ge\beta_0(\alpha)$. Hence the `diagonal' net $I(F_{\gamma})\to 
I(X)$ $(\tau_n)$. 
On the other hand, if $f\in I(X)$ and $\epsilon>0$ is given then a simple 
topological argument shows that 
there exists $\alpha_0$ such that for $\alpha\ge\alpha_0$, 
$N_{\alpha}\subseteq\{x\in \Gamma(A):\vert f(x)\vert<\epsilon\}$. Thus for 
$\alpha\ge\alpha_0$, $$\Vert f+I(F_{\beta(\alpha)})\Vert <k\epsilon$$ for all 
$\beta(\alpha)$, by the Helson property (with constant $k$). Hence 
$I(F_{\gamma})\to I(X)$ $(\tau_u)$, 
using [\the\Id; 2.1], 
and so 
$$I(F_{\gamma})\to I(X)~~(\tau_r).$$ 
Since $I\subseteq L(X)$, 
and $L(X)$  is a strict subset of $I(X)$, 
we have $I\ne I(X)$, so $\tau_r$ is not 
Hausdorff. Q.E.D. 
\bigskip 
\noindent Recall that for a Banach algebra, $\tau_r$ is Hausdorff if and only 
if $\tau_n$ has the normality property [\the\Id; 2.12]. Thus Theorem 1.9 can 
be rephrased as saying that for semisimple commutative Banach algebras with 
the Helson property, spectral synthesis holds if and only if $\tau_n$ has the 
normality property. In this form, Theorem 1.9 is evidently closely related to 
Beckhoff's result [\the\Beck; Proposition 3, Theorem 6] that for semisimple 
commutative Banach algebras with the `distance property', spectral synthesis 
holds if and only if $\tau_w$ has the normality property. 
\bigskip 
\noindent We now need to check that the change of definition of spectral 
synthesis has not spoiled the results of the previous paper [\the\Syn]. Recall 
the definitions of weak and strong spectral synthesis from [\the\Syn], given 
at the beginning of this section. We extend these to the non-unital case 
simply by dropping the requirement for $A$ to be unital. Note that the Baire 
property, axiom (ii), is superfluous for weak spectral synthesis, since 
$Prime(A)$ is always sober, and hence automatically a Baire space when it is 
locally compact. 

First we observe that the change in definition of spectral synthesis allows us 
to remove the separability restriction from [\the\Syn; 1.2]. 
\bigskip 
\noindent{\bf Proposition 1.10} {\sl Let $A$ be a Banach algebra. If $A$ has 
strong spectral synthesis then $A$ has spectral synthesis. If $A$ has spectral 
synthesis then $A$ has weak spectral synthesis.} 
\bigskip 
\noindent{\bf Proof.} Suppose that $A$ has strong spectral synthesis. Then 
every proper closed ideal of $A$ is strongly semisimple, so 
$Max^s(A)=Prim^s(A)$ by Lemma 1.1. Thus condition (i)$'$ of Definition 1.2 
holds, since $Max^s(A)$ is locally compact by Remark (a) after Definition 1.2. 
Condition (ii)$'$ follows by an easy modification of [\the\Syn; 1.1], and 
condition (iii)$'$ is trivial. Thus $A$ has spectral synthesis. 

Now suppose that $A$ has spectral synthesis. Then every closed ideal of $A$ is 
semisimple, so $Prime(A)=Prim^s(A)$ by Lemma 1.1. The remarks after Definition 
1.2 now establish that $A$ has weak spectral synthesis. Q.E.D. 
\bigskip 
\noindent Now that weak spectral synthesis has been defined for non-unital 
Banach algebras, we can observe that [\the\Syn; Theorem 2.9] extends to the 
non-unital case. The only place in the proof where the identity was used was 
in the appeal to [\the\Syn; Corollary 2.8], but that corollary is valid for 
non-unital Banach algebras (simply adjoin an identity). Thus we have the 
following. 
\bigskip 
\noindent {\bf Theorem 1.11} {\sl Let $A$ be a Banach algebra and suppose that 
$\tau_{\infty}$ is Hausdorff on $Id(A)$. Then $A$ has weak spectral 
synthesis.} 
\bigskip 
\noindent We do not know in general whether weak spectral synthesis implies 
spectral synthesis. It was remarked in [\the\Syn] that an 
infinite-dimensional, simple, radical Banach algebra (if there is one) would 
have weak spectral synthesis, but not spectral synthesis. Thus even in the 
commutative case it could be that weak spectral synthesis does not imply 
spectral synthesis. 

The next result, motivated by discussion with Peter Dixon, shows that the 
general problem of showing that $\tau_{\infty}$ being Hausdorff on $Id(A)$ 
implies that $A$ has spectral synthesis is more difficult than might appear. 
One would first have to show that there are no algebraically simple, radical 
Banach algebras---a famous old problem. Such an algebra would not have 
spectral synthesis, but would have $\tau_{\infty}$ Hausdorff, as we now show. 
The following lemma is certainly not new. 
\bigskip 
\noindent {\bf Lemma 1.12} {\sl Let $A$ be an algebraically simple Banach 
algebra of dimension greater than one, and let $x$ be a non-zero element of 
$A$. Set $$AxA=\left\{\sum_{i=1}^n a_ixb_i: a_1,\ldots , a_n, b_1,\ldots , 
b_n\in A\right\}.$$ Then $AxA=A$.} 
\bigskip 
\noindent {\bf Proof.} The set $AxA$ is an ideal. If $AxA=\{ 0\}$ for some 
non-zero $x$, then the ideal $I=\{ x\in A:AxA=\{0\}\}$ would be non-zero. 
Hence $I$ would equal $A$. This would imply that $A^3=\{0\}$, and this would 
lead immediately to the existence of non-zero proper ideals. Q.E.D. 
\bigskip 
\noindent {\bf Theorem 1.13} {\sl Let $A$ be an algebraically simple Banach 
algebra of dimension greater than one. Then the topology $\tau_{\infty}$ is 
Hausdorff on $Id(A)$.} 
\bigskip 
\noindent {\bf Proof.} The point $A$ is always $\tau_{\infty}$-closed in 
$Id(A)$, so we need only show that $\{0\}$ is $\tau_{\infty}$-closed. To do 
this we shall show that for each non-zero $x\in A$ there is a constant $C>0$ 
such that $\rho(x)\ge C/k^3$ for all non-zero $\rho\in S_k$. Set $I=AxA$. Then 
$I=A$ by Lemma 1.12, so $I^2=A$ (since $A^2\ne 0$). Thus there exist $a_i, 
b_i, c_i\in A$ $(1\le i\le n)$ such that $$x=\sum_{i=1}^n a_ixb_i xc_i .$$ Let 
$\rho\in S_k$. Then $$0<\rho(x)\le\sum_{i=1}^n 
\rho(a_i)\rho(x)\rho(b_i)\rho(x)\rho(c_i) 
\le \rho(x)^2\sum_{i=1}^n k^3\Vert a_i\Vert\, \Vert b_i\Vert\, \Vert 
c_i\Vert.$$ Thus $$1\le \rho(x) \sum_{i=1}^n k^3\Vert a_i\Vert\, \Vert 
b_i\Vert\, \Vert c_i\Vert$$ so $$C=\left\{\sum_{i=1}^n \Vert a_i\Vert\, \Vert 
b_i\Vert\, \Vert c_i\Vert\right\}^{-1}$$ is the required constant. Q.E.D. 
\bigskip 
\noindent Although we do not know in general whether weak spectral synthesis 
implies spectral synthesis, we shall show in a moment (Theorem 1.17) that this 
implication holds for semisimple PI-Banach algebras, and hence in particular 
for semisimple commutative Banach algebras. First we need to show that the 
various forms of spectral synthesis have reasonable stability properties. 
\bigskip 
\noindent {\bf Proposition 1.14} {\sl Let $A$ be a Banach algebra with weak 
spectral synthesis. If $I\in Id(A)$ and $J\in Id(I)$ then $J\in Id(A)$.} 
\bigskip 
\noindent {\bf Proof.} We show that $J$ is a union of ideals of $A$. 

For $a\in A$, let $I_a$ be the smallest closed, two-sided ideal of $A$ 
containing $a$. Evidently $I_a\subseteq I$ if $a\in I$. By condition (iv) of 
the definition of weak spectral synthesis, $I_a=\bigcap\{ P\in Prime(A):a\in 
P\}$. It is clear that $I_a\supseteq (AaA)^-$. On the other hand, if $P\in 
Prime(A)$ and $a\notin P$ then $AaA\not\subseteq P$. Thus $P\supseteq 
AaA\Leftrightarrow P\supseteq I_a$, so condition (iv) implies that 
$(AaA)^-=I_a$. 

Note too that $(K^3)^-=K$ for all $K\in Id(A)$, since a closed prime ideal 
contains $(K^3)^-$ if and only if it contains $K$. Hence 
$$I_a=((I_a)^3)^-=(I_a(AaA)^- I_a)^-=(I_a(AaA)I_a)^-$$ 
$$=((I_aA)a(AI_a))^- 
\subseteq (I_aaI_a)^-\subseteq (AaA)^-=I_a.$$ 
Thus $I_a=(I_aaI_a)^-$. It follows that 
$$J\supseteq\bigcup_{a\in J}(IaI)^-\supseteq\bigcup_{a\in J}(I_aaI_a)^-$$ 
$$=\bigcup_{a\in J}I_a\supseteq\bigcup_{a\in J}\{ a\}=J.$$ 
Hence $J=\bigcup_{a\in J}I_a$, so $J\in Id(A)$. Q.E.D. 
\bigskip 
\noindent {\bf Corollary 1.15} {\sl Let $A$ be a Banach algebra with spectral 
synthesis (or weak, or strong spectral synthesis) and let $J\in Id(A)$. Then 
$J$ and $A/J$ have spectral synthesis (or weak, or strong spectral synthesis, 
respectively).} 
\bigskip 
\noindent {\bf Proof.} First we give the proof for spectral synthesis. 
Conditions (i)$'$ and (iii)$'$ of Definition 1.2 follow for $A/J$ from the 
canonical homeomorphisms between $Prim(A/J)$ and $\{ P\in Prim(A):P\supseteq 
J\}$, and $Id(A/J)$ and $\{I\in Id(A):I\supseteq J\}$, respectively [\the\Ri; 
2.6.6]. Condition (ii)$'$ follows from the fact that the first of these 
homeomorphisms preserves the quotient norms. 

The spectral synthesis for $J$ is slightly less obvious. Condition (i)$'$ for 
$J$ follows from the homeomorphism $P\mapsto P\cap J$ between $\{ P\in 
Prim(A):P\not\supseteq J\}$ and $Prim(J)$. Condition (ii)$'$ for $J$ follows 
by using this homeomorphism and condition (ii)$'$ for $A$, and noting that 
$$\Vert b+P\Vert\le\Vert b+P\cap J\Vert\ \ \ \ \ \ \ \ (b\in J,\, P\in 
Prim(A)).$$ Condition (iii)$'$ for $J$ follows from condition (iii)$'$ for $A$ 
together with Proposition 1.14. 

The proof for weak spectral synthesis follows in the same way. 
For strong spectral synthesis, the only additional feature is the need to 
establish condition (ii), the Baire property for closed subsets. Each closed 
subset of $Max(A/J)$ is homeomorphic to a closed subset of $Max(A)$, and hence 
has the Baire property by condition (ii) for $A$. Now, suppose that $M$ is a 
closed subset of $Max(J)$, and that $(X_i)_{i\ge 1}$ is a countable collection 
of dense open subsets of $M$. Let $N$ be the closure of $M$ in $Max(A)$ 
(regarding $Max(J)$ as a subset of $Max(A)$ in the usual way). Then each $X_i$ 
is dense in $N$, and open since $Max(J)$ is open in $Max(A)$. But $N$ has the 
Baire property by condition (ii) for $A$, so $\bigcap_{i=1}^{\infty}X_i$ is 
dense in $N$, and thus in $M$. It follows that $M$ has the Baire property. 
Q.E.D. 
\bigskip 
\noindent It would be interesting to know whether the extension property holds 
for spectral synthesis, that is, suppose that $A$ is a Banach algebra and that 
$J\in Id(A)$ with $J$ and $A/J$ both having spectral synthesis. Does $A$ have 
spectral synthesis? We do not know the answer even in the special cases when 
$A$ is commutative, or when $A$ is the unitization of $J$. The best that we 
have is the following. 
\bigskip 
\noindent{\bf Proposition 1.16} {\sl Let $A$ be a Banach algebra with $J\in 
Id(A)$, and suppose that $J$ and $A/J$ have spectral synthesis. Then every 
closed ideal of $A$ is semisimple.} 
\bigskip 
\noindent{\bf Proof.} We introduce the following temporary notation. For a 
Banach algebra $B$, and $I\in Id(B)$, let $h_B(I)=\{P\in Prim(B):P\supseteq 
I\}$. Now let $\pi:A\to A/J$ denote the quotient map. Let $I,\, K\in Id(A)$ 
with $h_A(I)=h_A(K)$. We shall show that $I=K$. First note that the assumption 
on $I$ and $K$ implies that $h_{A/J}(\pi(I))=h_{A/J}(\pi(K))$ and hence that 
$\pi(I)=\pi(K)$, by spectral synthesis for $A/J$. On the other hand, we also 
have that $h_J(I\cap J)=h_J (K\cap J)$, and hence that $I\cap J=K\cap J$, by 
spectral synthesis for $J$. 

Now let $b\in I$. Then there exists $c\in K$ such that $\pi(b)=\pi(c)$, and 
hence such that $b-c\in J$. But this means (see the proof of Proposition 1.14) 
that $$(b-c)\in (J(b-c)J)^-\subseteq (JbJ-JcJ)^-$$ $$\subseteq (JIJ -JKJ)^-$$ 
$$\subseteq (I\cap J-K\cap J)^-=K\cap J.$$ Hence $b=c-(c-b)\in K+K\cap J=K$. 
Thus $I\subseteq K$. Similarly, $K\subseteq I$. Q.E.D. 
\bigskip 
\noindent Note that the proof of Proposition 1.16 only uses the fact that 
every closed ideal in $J$ and $A/J$ is semisimple. Proposition 1.16 raises the 
interesting question of whether a Banach algebra for which every closed ideal 
is semisimple automatically has spectral synthesis. We shall see that this is 
the case for the class of examples which we consider in the next section. In 
general the answer to the question is not known even for commutative Banach 
algebras. 

We can now show that for semisimple PI-Banach algebras, weak spectral 
synthesis is equivalent to spectral synthesis. Let us introduce some notation. 
If $A$ is a Banach algebra, and $P$ is the kernel of an irreducible 
representation of $A$ of dimension $k$ then $P$ has codimension $k^2$. Let 
$X_k$ $(k\ge 1)$ be the set of primitive ideals of $A$ of codimension less 
than or equal to $k^2$, and let $X_0$ be the empty set. Kaplansky showed that 
$X_k$ is a $\tau_w$-closed subset of $Prim(A)$ [\the\Kap; p.237]. 

Recall that an algebra is said to be a {\sl PI-algebra} if there is a 
non-trivial polynomial identity satisfied by the elements of $A$ (for example, 
all $X,Y,Z\in M_2({\bf C})$, the $2\times 2$ complex matrices, satisfy 
$[[X,Y]^2,Z]=0$, where $[\,.\, ,\, .\, ]$ denotes the additive commutator). We 
shall say that $A$ is a PI-algebra of {\sl degree m} if $A$ satisfies a 
polynomial of degree $m$, but of no smaller degree. If $A$ is a Banach algebra 
satisfying a polynomial identity of degree $m$ then every primitive ideal of 
$A$ is maximal and every irreducible representation of $A$ has dimension less 
than or equal to $m/2$, see [\the\Cohn; 10.4.6] and [\the\Pal; 7.1.16]. If $A$ 
is a semisimple PI-Banach algebra of degree $m$ then $A$ has an irreducible 
representation of degree $m/2$. 
\bigskip 
\noindent{\bf Theorem 1.17} {\sl Let $A$ be a PI-Banach algebra with weak 
spectral synthesis. If $A$ is semisimple then $A$ has spectral synthesis.} 
\bigskip 
\noindent{\bf Proof.} It is enough to show that every proper closed prime 
ideal of $A$ is primitive. 

Suppose first that all the primitive ideals of $A$ are of the same finite 
codimension $n^2$. Let $R$ be any { proper} 
closed ideal of $A$ of codimension not 
greater than $n^2$. Then $A/R$ is finite-dimensional, and semiprime, so $A/R$ 
is semisimple by Wedderburn's theorem [\the\Pal; 8.1.1]. Thus 
$R\in Prim(A)$, 
by the assumption on dimensions. 
{ This shows that $Prim(A)$ coincides with the set of proper
closed ideals of $A$ of 
codimension not greater than $n^2$.}
{ By [\the\Id; 5.1], this set is $\tau_n$-closed in $Id'(A)$ (the set 
of all proper, closed, two-sided ideals of $A$).
Also,  since $\tau_w$ has the normality 
property on $Prime(A)$, $\tau_w$ and $\tau_n$ coincide on $Prime(A)$.
Thus $Prim(A)$ is $\tau_w$-closed in $Prime(A)$.}
On the other hand, $Prim(A)$ is 
$\tau_w$-dense in $Prime(A)$ because $A$ is semisimple [\the\Id; 1.1]. Thus 
$Prim(A)=Prime(A)$. 

Now suppose that $A$ is any PI-Banach algebra of degree $2n$, with weak 
spectral synthesis. Let $J_i=\bigcap X_i$ $(0\le i\le n)$ (so $J_n=\{0\}$, and 
$J_0=A$). Let $P\in Prime(A)$ and let $m\in\{0,\ldots, n-1\}$ be the smallest 
number such that $P\not\supseteq J_m$. Thus $P\supseteq J_{m+1}$. Let 
$\pi:A\to A/J_{m+1}$ be the quotient map. Set $B=A/J$, $Q=\pi(P)$, and 
$I=\pi(J_m)$. Then $Q\in Prime(B)$ and $Q\in Prim(B)$ if and only if $P\in 
Prim(A)$. Since $Q\not\supseteq I$, $Q\cap I\in Prime(I)$, and $Q\cap I\in 
Prim(I)$ if and only if $Q\in Prim(B)$. But $Prim(I)\cong X_m\setminus 
X_{m-1}$, and all the primitive ideals of $I$ are of the same finite 
codimension $m^2$. Furthermore $I$ has weak spectral synthesis by Corollary 
1.15. Thus $Prime(I)=Prim(I)$ by the previous paragraph, so $Q\cap I\in 
Prim(I)$. Hence $P\in Prim(A)$, as required. Q.E.D. 

\bigskip 
\noindent Combining Theorems 1.11 and 1.17, we get the following. The unital 
case was established in [\the\Syn; 2.10]. 
\bigskip 
\noindent{\bf Corollary 1.18} {\sl Let $A$ be a PI-Banach algebra with 
$\tau_{\infty}$ Hausdorff on $Id(A)$. Then $A$ has spectral synthesis.} 
\bigskip 
\noindent{\bf Proof.} By Theorem 1.11 and Theorem 1.17 we need only show that 
$A$ is semisimple. But $A$ is semiprime, by Theorem 1.11, so if the Jacobson 
radical of $A$ were non-zero then there would be a non-zero quasi-nilpotent 
element in the centre of $A$ by Rowen's theorem [\the\Cohn; 10.7.5] (every 
non-zero ideal in a semiprime PI-algebra has non-zero intersection with the 
centre). As was observed in Remark (ii) after Proposition 2.7 of [\the\Syn], 
such an element is all that is needed to construct a net of norms in $S_1(A)$ 
converging to a seminorm with a non-zero kernel, and this would contradict the 
Hausdorffness of $\tau_{\infty}$. Thus $A$ is semisimple. Q.E.D. 
\bigskip 
\noindent A partial converse to Corollary 1.18 was established in [\the\Syn; 
3.10]: if $A$ is a separable, unital PI-Banach algebra with spectral synthesis 
then $\tau_{\infty}$ is Hausdorff. 
\bigskip 
\bigskip 
\noindent {\bf 2. Haagerup tensor products of C$^*$-algebras} 
\bigskip 
\noindent In this section we start looking at spectral synthesis and the 
topology $\tau_{\infty}$ for the Haagerup tensor product $A\otimes_h B$ of two 
C$^*$-algebras $A$ and $B$. We show that if $A$ and $B$ are separable then 
every closed prime ideal of $A\otimes_hB$ is primitive, just as for 
C$^*$-algebras. For general $A$ and $B$ we show that $A\otimes_h B$ has 
spectral synthesis if and only if $\tau_{\infty}$ is Hausdorff on 
$Id(A\otimes_h B)$. 

For C$^*$-algebras $A$ and $B$, the {\sl Haagerup norm} $\Vert\, .\, \Vert_h$ 
is defined on an element $x$ in the algebraic tensor product $A\otimes B$ by 
$$\Vert x\Vert_h=\ \inf\ 
\left\Vert\sum_{i=1}^na_ia_i^*\right\Vert^{1/2}\left\Vert\sum_{i=1}^nb^*_ib_i 
\right\Vert^{1/2},$$ 
where the infimum is taken over all possible representations of $x$ as a 
finite sum $x=\sum_{i=1}^n\ a_i\otimes b_i$, $a_i\in A$, $b_i\in A$. The 
completion of $A\otimes B$ in this norm is called the {\sl Haagerup tensor 
product}. The norm is a Banach algebra norm; introductory information on the 
Banach algebra $A\otimes_h B$ can be found in [\the\Ble], [\the\ASS], and 
[\the\AKSS]. Since this is still something of a specialist area, we assemble a 
few useful facts. The first two are the most remarkable. 
\medskip 
(i) Injectivity. If $A$ and $B$ are C$^*$-algebras, and $C$ and $D$ are 
C$^*$-subalgebras of $A$ and $B$ respectively, then the algebraic embedding 
$C\otimes D\to A\otimes_h B$ extends to an isometric embedding of $C\otimes_h 
D$ in $A\otimes_h B$ [\the\PS]. 

(ii) Exactness. If $A$ and $B$ are C$^*$-algebras and $I\in Id(A)$, $J\in 
Id(B)$, then there is a short exact sequence 
$$0\to I\otimes_h B+A\otimes_h J\to A\otimes_h B\to A/I\otimes_h B/J\to 0.$$ 

(iii) If $A$ and $B$ are both infinite-dimensional C$^*$-algebras then 
$A\otimes_h B$ is not bicontinuously isomorphic to a Banach subalgebra of the 
algebra of bounded operators on a Hilbert space [\the\Blec; 2.2]. 

(iv) The natural involution on the algebraic tensor product $A\otimes B$ of 
two C$^*$-algebras $A$ and $B$ does not in general extend to an involution on 
$A\otimes_h B$, but if $A$ and $B$ are abelian then the involution does 
extend, and $A\otimes_h B$ is a Banach $^*$-algebra. 

(v) The Banach algebra $A\otimes_h B$ is semisimple [\the\ASS; 5.16]. 

(vi) The Haagerup tensor product $A\otimes_h B$ of two abelian C$^*$-algebras 
$A$ and $B$ is bicontinuously isomorphic to their projective tensor product 
$A\hat\otimes B$, see [\the\ASS; p.112]. Thus $A\otimes_h B$ is in this case a 
completely regular, semisimple commutative Banach $^*$-algebra. 
\medskip 
\noindent We shall use these various properties without further comment. 
\bigskip 
The first thing to do now is to consider the various different spaces of 
`primitive' ideals of $A\otimes_h B$. Recall that $Prim^s(A\otimes_h B)$ is 
the sobrification of $Prim(A\otimes_h B)$ (the space of primitive ideals of 
$A\otimes_h B$) and consists of the semisimple prime ideals of $A\otimes_h B$. 
It follows from [\the\ASS; 5.10, 2.6, and 5.16] that every closed prime ideal 
of $A\otimes_h B$ is semisimple. Thus $Prim^s(A\otimes_h B)=Prime(A\otimes_h 
B)$, the set of proper, closed, prime ideals of $A\otimes_h B$, and hence the 
notions of spectral synthesis and weak spectral synthesis coincide for the 
class of Haagerup tensor products of C$^*$-algebras. Furthermore 
$Prim^s(A\otimes_h B)$ is homeomorphic to $Prime(A)\times Prime(B)$, by 
[\the\ASS; 5.10] and [\the\AKSS; 1.5], so $Prim^s(A\otimes_h B)$ is always 
locally compact. Thus condition (i)$'$ of spectral synthesis is always 
satisfied for $A\otimes_h B$. The same is true for condition (ii)$'$, by 
[\the\AKSS; 3.4]. Thus the only condition that needs to be checked for 
spectral synthesis is (iii)$'$ --- whether every closed ideal of $A\otimes_h 
B$ is semisimple. 

Notice that we can therefore solve the extension problem for algebras in this 
class with spectral synthesis. It follows from the remark after Proposition 
1.16 that if $M$ is an algebra in this class with a closed ideal $J$ such that 
both $J$ and $M/J$ have spectral synthesis then $M$ has spectral synthesis. 

The other space of `primitive' ideals that we have to bear in mind is 
$Prim_*(A\otimes_h B)$, defined as follows. A representation $\pi$ of 
$A\otimes_h B$ on a Hilbert space $H$ is a {\sl $*$-representation} if 
$\pi(a^*\otimes b^*)=\pi(a\otimes b)^*$ for $a\in A,\ b\in B$. Let 
$Prim_*(A\otimes_h B)$ be the space of kernels of topologically irreducible 
$*$-representations of $A\otimes_h B$, with the hull-kernel topology. Elements 
of $Prim_*(A\otimes_h B)$ are prime and semisimple [\the\ASS; 5.13(ii), 5.10, 
5.16], so $Prim_*(A\otimes_h B)\subseteq Prim^s(A\otimes_h B)$, with equality 
holding if $A$ and $B$ are separable [\the\ASS; 5.15]. In Corollary 2.3 we 
shall show that, for separable $A$ and $B$, these two spaces also coincide 
with $Prim(A\otimes_h B)$. Caution is required in using results from 
[\the\ASS] and [\the\AKSS] because the definition of `primitive' used there 
corresponds to our definition of elements of $Prim_*(A\otimes_h B)$ (in 
particular, the statement of [\the\ASS; 5.13(i)] is identical to that of our 
Theorem 2.2!). 

Ideals of $A\otimes_h B$ of the form $ I\otimes_h J$ $(I\in Id(A),\ J\in Id 
(B))$ are called {\sl product ideals}. Clearly each product ideal is generated 
by the simple tensors (i.e. tensors of the form $a\otimes b$, $a\in A,\ b\in 
B$) which it contains. An ideal which is generated by the product ideals which 
it contains is called a {\sl lower} ideal, see [\the\ASS; 6.7]. To define 
upper ideals, we use the observation of Blecher's [\the\Ble] that if $A$ and 
$B$ are two C$^*$-algebras then the Haagerup tensor product $A\otimes_h B$ is 
canonically isomorphic to a dense subalgebra of the minimal C$^*$-tensor 
product $A\otimes_{min}B$. A closed ideal in $A\otimes_h B$ of the form $J\cap 
A\otimes_h B$ $(J\in Id(A\otimes_{min} B))$ is called an {\sl upper ideal} 
[\the\ASS; \S6]. It was shown in [\the\ASS; 6.7] that an ideal in 
$Id(A\otimes_h B)$ is an upper ideal if and only if it is an intersection of 
elements of $Prim_*(A\otimes_h B)$. Each closed ideal $I$ is contained in a 
unique smallest upper ideal $I^u$, and contains a unique largest lower ideal 
$I_l$. {\sl Spectral synthesis} is said to hold at $I$ if $I^u=I_l$ [\the\ASS; 
p.140], and $A\otimes_h B$ was said in [\the\ASS] to have spectral synthesis 
if every closed ideal has spectral synthesis. 
It is known that every product ideal has spectral synthesis; indeed the sum of 
a finite number of product ideals is itself a closed ideal [\the\ASS; 3.8] and 
has spectral synthesis [\the\ASS; 6.4]. Since every primitive ideal is closed 
and prime, and every closed prime ideal is a sum of two product ideals 
[\the\ASS; 5.9], it follows that every primitive ideal has spectral synthesis. 
In particular, every primitive ideal of $A\otimes_h B$ is an upper ideal. On 
the other hand, since every element of $Prim_*(A\otimes_h B)$ is semisimple, 
it follows that an ideal of $A\otimes_h B$ is an upper ideal if and only if it 
is semisimple. 

Putting all this together, we have established the following. 
\bigskip 
\noindent {\bf Proposition 2.1} {\sl Let $A$ and $B$ be C$^*$-algebras. Then 
the Haagerup tensor product $A\otimes_h B$ has spectral synthesis in the sense 
of this paper if and only every closed ideal of $A\otimes_h B$ has spectral 
synthesis in the sense of {\rm [\the\ASS]}.} 
\bigskip 
\noindent {\bf Examples} (i) If $A$ and $B$ are C$^*$-algebras with $Id(A)$ or 
$Id(B)$ finite then $A\otimes_h B$ has spectral synthesis, by [\the\ASS; 5.3, 
6.4]. Thus there are lots of non-commutative Banach algebras which are not 
C$^*$-algebras but which have spectral synthesis, e.g. take the Haagerup 
tensor product of any infinite-dimensional C$^*$-algebra with the compact 
operators, or with the Fermion algebra. 

(ii) Let $X$ and $Y$ be countable, locally compact, Hausdorff spaces. Then 
$C_0(X)\otimes_h C_0(Y)$ has spectral synthesis. To see this, we argue as 
follows. First note that $C_0(X)\otimes_h C_0(Y)$ is a closed ideal in 
$C(X')\otimes_h C(Y')$ (where $X'$ and $Y'$ are the one-point 
compactifications of $X$ and $Y$ respectively), so by Corollary 1.15 we may 
assume that $X$ and $Y$ are compact. It now follows by transfer methods, see 
[\the\GM; 11.1.3], that there is a locally compact abelian group $G$ such that 
$C(X)\otimes_h C(Y)$ is bicontinuously isomorphic to a quotient of $L^1(G)$ by 
a semisimple ideal $J$. Thus $X\times Y$ 
is homeomorphic to the hull of $J$ in $\widehat G$. But every countable closed 
subset of $\widehat G$ is a set of synthesis, see [\the\Rud; 7.2.5], from 
which it follows that $C(X)\otimes_h C(Y)$ has spectral synthesis. 

(iii) If $D$ is the Cantor set then the commutative Banach algebra 
$C(D)\otimes_h C(D)$ does not have spectral synthesis (because the projective 
tensor product $C(D)\hat\otimes C(D)$, with which it is isomorphic, does not 
have spectral synthesis, see [\the\GM; 11.2.1]). 

\bigskip 
We now show that if $A$ and $B$ are separable C$^*$-algebras then every closed 
prime ideal of $A\otimes_hB$ is primitive. For a C$^*$-algebra $A$, let 
$\tilde A$ denote the unitization of $A$ (that is, $A$ itself if $A$ is 
unital, and $A$ with an identity adjoined otherwise). 
\bigskip 
\noindent {\bf Theorem 2.2} {\sl Let $A$ and $B$ be C$^*$-algebras, and let 
$P\in Prim(A)$ and $Q\in Prim(B)$. Then $P\otimes_h B+A\otimes_h Q\in 
Prim(A\otimes_hB)$.} 
\bigskip 
\noindent{\bf Proof.} Since $A\otimes_h B/(P\otimes_h B+A\otimes_h Q)\cong 
A/P\otimes_h B/Q$ [\the\ASS; 2.5], it is enough to show that $A\otimes_h B$ is 
a primitive Banach algebra under the hypothesis that $A$ and $B$ are primitive 
C$^*$-algebras. Furthermore, since $A\otimes_h B$ is a closed two-sided ideal 
of $\tilde A \otimes_h \tilde B$, it is enough to show that $\tilde A 
\otimes_h \tilde B$ is primitive. Thus we may assume, at the outset, that $A$ 
and $B$ are unital. 

Let $M$ and $N$ be maximal left ideals of $A$ and $B$ respectively such that 
the largest two-sided ideal of $A$ (respectively, $B$) contained in $M$ 
(respectively $N$) is the zero ideal. Then the left ideal $M\otimes_h 
B+A\otimes_h N$ is closed in $A\otimes_h B$ [\the\ASS; 2.6]. Let $S$ be any 
left ideal of $A\otimes_h B$ containing both $M\otimes_h B+A\otimes_h N$ and a 
non-zero closed, two-sided ideal of $A\otimes_h B$. We shall show that 
$S=A\otimes_h B$. 

Each non-zero closed two-sided ideal of $A\otimes_h B$ contains a non-zero 
simple tensor [\the\ASS; 4.6], so there is a non-zero simple tensor $a\otimes 
b\in S$ such that the closed two-sided ideal $I$ generated by $a\otimes b$ is 
contained in $S$. On the other hand, since $A$ and $B$ are primitive, 
$I_a\not\subseteq M$ and $I_b\not\subseteq N$, where $I_a$ and $I_b$ are the 
smallest closed, two-sided ideals of $A$ and $B$ containing $a$ and $b$ 
respectively. Thus there exist $c,d\in A$ and $e,f\in B$ such that $cad\notin 
M$ and $ebf\notin N$. 

Since $M$ and $N$ are maximal left ideals in $A$ and $B$ respectively there 
exist $g\in A$ and $m\in M$ such that $g(cad)+m=1_A$, and $h\in B$ and $n\in 
N$ such that $h(ebf)+n=1_B$. Hence 
$$\eqalign{(1_A\otimes 1_B)&=(gcad+m)\otimes (hebf+n)\cr 
&= gcad\otimes hebf+gcad\otimes n+m\otimes 1_B\cr 
&=(gc\otimes he)(a\otimes b)(d\otimes f)+(1_A\otimes n-m\otimes n)+m\otimes 
1_B\cr 
&\in I+(S-S)+S\subseteq S.\cr}$$ 
Hence $S=A\otimes_h B$ as required. Thus if $T$ is a maximal left ideal of 
$A\otimes_h B$ containing $M\otimes_h B+A\otimes_h N$ then the only two 
sided-ideal of $A\otimes_h B$ contained in $T$ is the zero ideal. Thus 
$A\otimes_h B$ is primitive. Q.E.D. 
\bigskip 
\noindent Theorem 2.2 leaves open the interesting question whether the 
topologically irreducible $*$-representation associated with the primitive 
ideal $P\otimes_h B+A\otimes_h Q$ is actually algebraically irreducible. 
\bigskip 
\noindent {\bf Corollary 2.3} {\sl Let $A$ and $B$ be separable 
C$^*$-algebras. Then $$Prim_*(A\otimes_h B)=Prim^s(A\otimes_h 
B)=Prim(A\otimes_h B).$$} 
{\bf Proof.} The equality of $Prim^s(A\otimes_h B)$ and $Prim_*(A\otimes_h 
B)$, for $A$ and $B$ separable, was established in [\the\ASS; 5.15]. If $R\in 
Prim(A\otimes_h B)$ then $R$ is, of course, prime and semisimple, so $R\in 
Prim^s(A\otimes_h B)$ (without any separability hypothesis). Conversely, 
suppose that $R\in Prim^s(A\otimes_h B)$. Then $R$ is closed and prime, so 
there exist closed, prime ideals $P\in Id(A)$ and $Q\in Id(B)$ such that 
$R=P\otimes_h B+A\otimes_h Q$ [\the\ASS; 5.9(iii)]. But closed, prime ideals 
are primitive in separable C$^*$-algebras, so $P\in Prim(A)$ and $Q\in 
Prim(B)$. Hence $R\in Prim(A\otimes_h B)$ by Theorem 2.2. Q.E.D. 
\bigskip 
\noindent We now show that spectral synthesis is equivalent to the 
Hausdorffness of $\tau_{\infty}$ for $A\otimes_h B$. We use the embedding of 
$A\otimes_h B$ in $A\otimes_{min} B$ mentioned earlier [\the\Ble]. With an eye 
to subsequent work we make the following general definition. 

Let $M$ be a Banach algebra, and suppose that there is a continuous norm 
$\gamma$ on $M$ such that $N$, the $\gamma$-completion of $M$, is a 
C$^*$-algebra. In this paper $M$ will of course be $A\otimes_h B$ and $N$ will 
be $A\otimes_{min} B$. Extending the earlier definition, we refer to those 
closed ideals in $M$ of the form $J\cap M$ $(J\in Id(N))$ as {\sl upper 
ideals}, and the set of upper ideals is denoted $Id^u(M)$. Note that if $I$ is 
an upper ideal of $M$ then in fact $I=J\cap M$ where $J$ is the closure of $I$ 
in $N$. 

\bigskip 
\noindent {\bf Definition 2.4} Let $M$ be a Banach algebra. We shall say that 
$M$ has {\sl property (P)} if $M$ satisfies the following conditions: 

(a) there is a continuous norm $\gamma$ on $M$ such that $N$, the 
$\gamma$-completion of $M$, is a C$^*$-algebra; 

(b) every primitive ideal of $M$ is an upper ideal, i.e. $Prim(M)\subseteq 
Id^u(M)$; 

(c) there is a subset $R$ of $M\cap N_{sa}$ such that each $a\in R$ is 
contained in a completely regular, commutative Banach $^*$-subalgebra $M_a$ of 
$M$ (where the norm and the involution on $M_a$ are those induced by $N$), and 
such that if $I\in Id^u(M)$ and $J\in Id(M)$ with $J\not\subseteq I$ then 
there exists $a\in R$ such that $M_a\cap J\not\subseteq I$. 
\bigskip 
\noindent Notice for later that if $M$ is a Banach algebra with property (P) 
and $I$ is an upper ideal of $M$ then $M/I$ also has property (P). 

In Theorem 2.7 we show that the Haagerup tensor product of two C$^*$-algebras 
has property (P). 
\medskip 
For the next lemma, recall Rickart's theorem, quoted before Theorem 1.8, which 
states that for any norm $\sigma$ on a completely regular, semisimple, 
commutative Banach algebra $A$, $\sigma(a)\ge\rho(a)$ for all $a\in A$, where 
$\rho$ denotes the uniform norm on $A$ [\the\Rick]. 
\bigskip 
\noindent{\bf Lemma 2.5} {\sl Suppose that $M$ is a Banach algebra with 
property (P) and that $a$ belongs to the subset $R$. Then 

(i) if $J\in Id^u(M)$ then $M_a\cap J$ is a semisimple ideal in $M_a$, 

(ii) $\gamma|_{M_a}$ is the uniform norm on $M_a$. } 
\bigskip 
\noindent{\bf Proof.} (i) Let $K$ be a closed ideal of $N$ such that $J=K\cap 
M$. Then $M_a/(M_a\cap J)$ is isomorphic to a commutative $^*$-subalgebra, $C$ 
say, of $N/K$. The completion of $C$ in $N/K$ is an abelian C$^*$-algebra, 
hence semisimple. Taking a completion reduces the spectral radius, so $C$ has 
no non-zero quasi-nilpotent elements. Thus $C$ is semisimple. 

(ii) Let $D\subseteq N$ be the $\gamma$-completion of $M_a$. Then $D$ is a 
commutative C$^*$-subalgebra of $N$. Let $\Delta(D)$ be the character space of 
$D$. For $\phi\in\Delta(D)$, let $\tilde\phi$ denote the character obtained by 
restricting $\phi$ to $M_a$. Then for $c\in M_a$, 
$$\gamma(c)=\sup\{\phi(c):\phi\in\Delta(D)\}=\sup\{\tilde\phi(c):\phi\in\Delta(D)\}\le 
\rho(c),$$ where $\rho$ denotes the spectral radius on $M_a$. But $M_a$ is 
semisimple by part (i) (taking $J=\{0\}$), so the spectral radius $\rho$ is 
the uniform norm on $M_a$. Hence $\rho(c)=\gamma(c)$, by Rickart's theorem, 
since $\gamma$ is also a norm. Q.E.D. 
\bigskip 
\noindent {\bf Proposition 2.6} {\sl Let $M$ be a Banach algebra with property 
(P) above. Then every upper ideal of $M$ is a $\tau_{\infty}$-closed point of 
$Id(M)$. If $M$ has spectral synthesis then $\tau_{\infty}$ is Hausdorff on 
$Id(M)$.} 
\bigskip 
\noindent {\bf Proof.} If $I$ is an upper ideal in $Id(M)$ then $M/I$ also has 
property (P), as we have observed, so we may assume that $I=\{0\}$. Next note 
that for any continuous norm $\rho$ on $M$ and for any $a\in R$, 
$\rho(a)\ge\gamma(a)$, by Rickart's theorem and Lemma 2.5(ii). Now suppose 
that $k\in{\bf N}$ and that $(\rho_{\alpha})$ is a net of norms in $S_k$, 
converging to a seminorm $\rho\in S_k$. Let $a\in R$ and let $b\in 
M_a\cap\ker\rho$. Then $\rho_{\alpha}(b)\to 0$. Hence $\gamma(b)=0$, so $b=0$. 
Thus $\ker\rho=\{0\}$, by condition (c) of property (P), so $\rho$ is a norm 
on $M$. Hence $\{ 0\}$ is $\tau_{\infty}$-closed, as required. 

Now suppose, additionally, that $M$ has spectral synthesis. Then every closed 
ideal of $M$ is semisimple, so condition (b) of property (P) implies that 
every closed ideal of $M$ is an upper ideal. Let $k\in{\bf N}$ and let 
$(\rho_{\alpha})$ be an increasing net of seminorms in $S_k$ with limit 
$\rho$. By [\the\Syn; 2.2] we must show that $J:=\sup\ker\rho_{\alpha}$ 
coincides with $\ker\rho$. It is always the case that $J\subseteq\ker\rho$. 
Let $a\in R$ and let $b\in M_a \cap \ker\rho$. For each $\alpha$, let 
$\sigma_{\alpha}$ be the uniform norm on the algebra $M_a/(M_a\cap 
\ker\rho_{\alpha})$, which is semisimple by Lemma 2.5(i), and let $\sigma$ be 
the uniform norm on the semisimple algebra $M_a/(M_a\cap J)$. Then by 
Rickart's theorem, $\rho_{\alpha}(b)\ge\sigma_{\alpha}(b)\ge\sigma(b)$ for 
each $\alpha$. Hence $\sigma(b)=0$, since $\rho_{\alpha}(b)\to 0$. Thus $b\in 
J$. It follows that $M_a\cap\ker\rho = M_a\cap J$, and hence that 
$J=\ker\rho$, by condition (c) of property (P). Q.E.D. 

\bigskip 
\noindent {\bf Theorem 2.7} {\sl Let $A$ and $B$ be C$^*$-algebras with 
Haagerup tensor product $A\otimes_hB$. Then every upper ideal is 
$\tau_{\infty}$-closed in $Id(A\otimes_h B)$, and $A\otimes_h B$ has spectral 
synthesis if and only if $\tau_{\infty}$ is Hausdorff on $Id(A\otimes_h B)$.} 
\bigskip 
\noindent {\bf Proof.} The fact that the Hausdorffness of $\tau_{\infty}$ 
implies spectral synthesis follows from Theorem 1.11 and the fact, already 
mentioned, that weak spectral synthesis and spectral synthesis are equivalent 
for $A\otimes_h B$. The rest of the theorem will follow from Proposition 2.6 
once we have established that $A\otimes_h B$ has property (P) of Definition 
2.4. 

We have already observed that conditions (a) and (b) of Property (P) are 
satisfied. Now let $R$ be the set $\{ a\otimes b: a\in A^+,\ b\in B^+\}$. Then 
for $a\otimes b\in R$ the commutative subalgebra $C^*(a)\otimes_h 
C^*(b)\subseteq A\otimes_h B$ is a completely regular Banach $^*$-algebra. 
Let $I\in Id^u(A\otimes_h B)$ and let $J\in Id(A\otimes_h B)$ with 
$J\not\subseteq I$. Then $J_l\not\subseteq I$ since $J\not\subseteq I$. But 
$J_l$ is generated by the product ideals which it contains [\the\ASS; 6.7(i)], 
and each product ideal is generated by its positive simple tensors, so there 
exists a simple tensor $a\otimes b\in J_l\setminus I$ with $a\in A^+$ and 
$b\in B^+$. This establishes condition (c) of property (P). 
Hence $A\otimes_h B$ has property (P), as required. Q.E.D. 
\bigskip 
\bigskip 
\noindent {\bf 3. The topology $\tau_r$ on $Id(A\otimes_h B)$} 
\bigskip 
\noindent In this section we consider the topology $\tau_r$ on $Id(A\otimes_h 
B)$, showing that if $Prime(A)$ and $Prime(B)$ are compact T$_1$-spaces, then 
$A\otimes_hB$ has spectral synthesis if and only if $\tau_r$ is Hausdorff on 
$Id(A\otimes_h B)$. The proof uses the theory of continuous lattices and 
bitopological spaces. 

We begin by looking at the lattice $Id^u(A\otimes_h B)$ of upper ideals. 
Recall that these are precisely the semisimple ideals of $A\otimes_h B$. Thus 
they are 
in bijective correspondence with the open (or closed) subsets of 
$Prim(A\otimes_h B)$, using the hull-kernel process. 

We need to recall the definition of a continuous lattice. Let $L$ be a 
complete lattice and let $x, y\in L$. Then $x$ is {\sl way-below} $y$, written 
$x\ll y$, if whenever $S$ is a subset of $L$ with $y\le\sup S$ there is a 
finite subset $F$ of $S$ such that $x\le \sup F$, see [\the\Comp]. If for each 
$y\in L$, $y=\sup\,\{ x\in L:x\ll y\}$ then $L$ is called a {\sl continuous 
lattice}. 
A continuous lattice carries three important topologies, as follows. The {\sl 
lower} topology on $L$ is the topology generated by the sets $\{ x\in 
L:x\not\ge y\}$ as $y$ varies through $L$ (in exact analogy with the lower 
topology on $Id(A)$ defined in the Introduction). A base for the {\sl Scott 
topology} on $L$ is given by sets of the form $\{ y\in L:y\gg x\}$ as $x$ 
varies through $L$. The {\sl Lawson} topology on $L$ is the join of the Scott 
and lower topologies. The Lawson topology on a continuous lattice is compact 
and Hausdorff [\the\Comp; III.1.10]. 

An element $p$ in a lattice $L$ is said to be {\sl prime} if whenever $x,y \in 
L$ with $x\wedge y\le p$ then either $x\le p$ or $y\le p$. Let $Prime(L)$ 
denote the set of proper prime elements of $L$ (in lattice theory it is 
conventional to include the `top' as well, but we are excluding it, in harmony 
with the convention in ring theory). { As in the earlier setting,
the lower topology coincides with the
hull-kernel topology on $Prime(L)$.}
A subset $S$ of $Prime(L)$ is said to be 
{\sl saturated} if $x\in S$ and $y\in Prime(L)$ with $y\le x$ implies that 
$y\in S$. 
We shall use the fact that a subset $S$ of $Prime(L)$ is open in the relative 
Scott topology if and only if $Prime(L)\backslash S$ is saturated and 
hull-kernel compact [\the\Comp; V.5.1]. For $x\in L$, let $hull(x)=\{p\in 
Prime(L): p\ge x\}$. The spectral theorem for a continuous, distributive 
lattice $L$ is that the correspondence $x\longleftrightarrow hull(x)$ is a 
lattice isomorphism between $L$ and the lattice of closed subsets of 
$Prime(L)$ in the lower topology [\the\Comp; V.5,5]. In particular, 
$x=\bigwedge hull(x)$ for each element $x$ in a continuous, distributive 
lattice. We shall also use the fact that if $x$ and $y$ are elements of a 
continuous distributive lattice $L$ then $x\ll y$ if and only if there exists 
a Scott open set $M$ such that $hull(x)\supseteq M\supseteq hull(y)$, see 
[\the\Comp; V.5.6(ii)]. The following result is standard, but it seems best to 
include the short proof. 
\bigskip 
\noindent{\bf Lemma 3.1} {\sl Let $L$ be a continuous distributive lattice. 
Let $x\in L$ and let $(x_{\alpha})$ be a net in $L$ converging to $x$ in the 
Scott topology. Let $V$ be a Scott open subset of $Prime(L)$ containing 
$hull(x)$. Then eventually $hull(x_{\alpha})\subseteq V$.} 
\bigskip 
\noindent {\bf Proof.} Suppose to the contrary. Then because 
$V^c:=Prime(L)\backslash V$ is hull-kernel compact, we may suppose, by passing 
to a subnet if necessary, that there exists $p\in V^c$ and a net 
$(p_{\alpha})$ with $p_{\alpha}\in V^c\cap hull(x_{\alpha})$ for each 
$\alpha$, such that $p_{\alpha}\to p$ as $\alpha\to\infty$, in the hull-kernel 
topology. Since $L$ is a continuous lattice, there exists $y\in L$ with $y\ll 
x$ and $y\not\le p$. The set $\{ z\in L:z\gg y\}$ is a Scott-open 
neighbourhood of $x$, so eventually $x_{\alpha}\gg y$. Hence eventually 
$x_{\alpha}\ge y$, so eventually $hull(x_\alpha)$ is contained in the set 
$hull(y)$ which is closed in the hull-kernel topology. Since $p\notin 
hull(y)$, we have a contradiction. Q.E.D. 
\bigskip 
We observed near the beginning of Section 2 that $Prim^s(A\otimes_h 
B)=Prime(A\otimes_h B)$ is homeomorphic to $Prime(A)\times Prime(B)$, and is 
therefore locally compact. On the other hand, the lattice $Id^u(A\otimes_h B)$ 
is isomorphic to the lattice of open subset of $Prim(A)$, which in turn is 
isomorphic to the lattice of open subsets of $Prim^s(A\otimes_h B)$. 
It follows from [\the\Comp; V.5.10], therefore, that $Id^u(A\otimes_h B)$ is a 
continuous distributive lattice. Furthermore, $Prime(Id^u(A\otimes_h B))$ is 
precisely $Prime(A\otimes_h B)$ [\the\ASS; 5.9]. The next proposition 
describes the lower, Scott, and Lawson topologies on $Id^u(A\otimes_h B)$. For 
$I\in Id^u(A\otimes_h B)$ and $x\in A\otimes_h B$, define $$\Vert 
x+I\Vert^u=\sup\{\Vert x+P\Vert: P\in hull(I)\},$$ where $hull(I)=\{ P\in 
Prime(A\otimes_h B): P\supseteq I\}$. We shall refer to the functions 
$I\mapsto\Vert x+I\Vert^u$ $(I\in Id^u(A\otimes_h B),\ x\in A\otimes_h B)$ as 
{\sl upper-norm functions}. Note that for $P\in Prime(A\otimes_h B)$, $\Vert 
x+P\Vert=\Vert x+P\Vert^u$ $(x\in A\otimes_h B)$. It was shown in [\the\AKSS; 
3.4] that for $x\in A\otimes_h B$, the function $P\mapsto\Vert x+P\Vert$ 
$(P\in Prime(A\otimes_h B))$ is lower semi-continuous on $Prime(A\otimes_h 
B)$. 
\medskip 
\noindent In Proposition 3.2 and Theorem 3.6 we shall use the hypothesis that 
$Prime(A)$ and $Prime(B)$ are compact, for the C$^*$-algebras $A$ and $B$. 
Note that for a C$^*$-algebra $A$, $Prime(A)$ is compact if and only if 
$Prim(A)$ is compact. This is easily seen from the `open-cover' version of 
compactness, using the natural extension of a cover of $Prim(A)$ to a cover of 
$Prime(A)$ by saturation, together with the fact that every proper closed 
prime ideal of $A$ is contained in a primitive ideal of $A$. One circumstance 
in which $Prim(A)$ is compact is, of course, when $A$ is unital. 
\bigskip 
\noindent{\bf Proposition 3.2} {\sl Let $A$ and $B$ be C$^*$-algebras with 
$Prime(A)$ and $Prime(B)$ compact in the hull-kernel topology. Let 
$C=A\otimes_hB$ be the Haagerup tensor product of $A$ and $B$ and let 
$Id^u(C)$ be the lattice of upper ideals of $C$. Then the lower, Scott, and 
Lawson topologies on $Id^u(C)$ are the weakest with respect to which the 
upper-norm functions $I\mapsto \Vert x+I\Vert^u$, $(x\in C, I\in Id^u(C))$, 
are respectively lower semi-continuous, upper semi-continuous, and 
continuous.} 
\bigskip 
\noindent{\bf Proof.} As temporary notation, let $\tau_{lc}$, $\tau_{uc}$, and 
$\tau_c$ be the weakest topologies on $Id^u(C)$ with respect to which the 
upper-norm functions are respectively lower semi-continuous, upper 
semi-continuous, and continuous. First we show that $\tau_{lc}$ coincides with 
the lower topology. 

For each $x\in C$, the set $\{ I\in Id^u(C): \Vert x+I\Vert^u>0\}$ is 
$\tau_{lc}$-open. Since these sets give a sub-base for the lower topology, it 
is immediate that $\tau_{lc}$ is stronger than the lower topology. Now let 
$(I_{\alpha})$ be a net in $Id^u(C)$ converging in the lower topology to some 
$I\in Id^u(C)$, and suppose, without loss of generality, that $x\in C$ with 
$\Vert x+I\Vert^u=1$. Let $1>\epsilon>0$ be given. Then there exists $P\in 
Prim(C)$ with $P\supseteq I$ such that $\Vert x+ P\Vert>1-{\epsilon\over 2}$. 
By the lower semicontinuity of the upper-norm functions on $Prim(C)$ 
[\the\AKSS; 3.4], there is a neighbourhood $X$ of $P$ in $Prim(C)$ such that 
$\Vert x+Q\Vert>1-\epsilon$ for all $Q\in X$. Let $y\in C\setminus P$ be such 
that the open subset $Y=\{Q\in Prim(C):y\notin Q\}$ of $Prim(C)$ is a 
neighbourhood of $P$ contained in $X$. Then $y\notin I$, so eventually 
$y\notin I_{\alpha}$. Thus eventually there is, for each $\alpha$, a primitive 
ideal $P_{\alpha}$ of $C$ containing $I_{\alpha}$ such that $y\notin 
P_{\alpha}$. Hence $P_{\alpha}\in Y\subseteq X$, so $\Vert 
x+P_{\alpha}\Vert>1-\epsilon$. Thus eventually $\Vert 
x+I_{\alpha}\Vert^u>1-\epsilon$. This shows that the lower topology is 
stronger than $\tau_{lc}$. Hence $\tau_{lc}$ is the lower topology on 
$Id^u(C)$. 

Next we show that $\tau_{uc}$ is weaker than the Scott topology on $Id^u(C)$. 
We begin by showing that the relative $\tau_{uc}$ topology is weaker than the 
relative Scott topology on $Prime(C)$. Recall that $Prime(C)$ is in bijective 
correspondence with $Prime(A)\times Prime(B)$, the correspondence being given 
by 
$$P\otimes_h B+ A\otimes_h B \longleftrightarrow (P,Q)\in Prime(A)\times 
Prime(B).$$ Recall too that the usual C$^*$-norm functions are upper 
semi-continuous on $Prime(A)$ and $Prime(B)$ in the relative Scott topologies 
[\the\HKMS; 7.2(b)]. Thus it follows from the remark after [\the\AKSS; 3.1] 
that the norm functions $I\mapsto \Vert x+P\Vert$, $(x\in C, P\in Prime(C))$, 
are upper semi-continuous on $Prime(C)$ when $Prime(C)$ is equipped with the 
product Scott topology from $Prime(A)\times Prime(B)$. We need only show, 
therefore, that the product Scott topology is weaker than the Scott topology 
on $Prime(C)$. This follows from the fact that if $X$ and $Y$ are compact, 
saturated subsets of $Prime(A)$ and $Prime(B)$ respectively then 
$$Z:=\{P\otimes_h B+A\otimes_h Q: (P,Q)\in (X\times Prime(B))\cup 
(Prime(A)\times Y)\}$$ 
is a compact saturated subset of $Prime(C)$. Hence the set 
$$Prime(C)\backslash Z\cong (Prime(A)\backslash X)\times (Prime(B)\backslash 
Y)$$ is open in $Prime(C)$ in the Scott topology. Since sets of the form 
$(Prime(A)\backslash X)\times (Prime(B)\backslash Y)$ form a base for the 
product Scott topology, it follows that the product Scott topology on 
$Prime(C)$ is weaker than the Scott topology. Hence the relative $\tau_{uc}$ 
topology is weaker than the relative Scott topology on $Prime(C)$. 

Now we extend this to the whole of $Id^u(C)$. Let $I\in Id^u(C)$ and let 
$(I_{\alpha})$ be a net in $Id^u(C)$ converging to $I$ in the Scott topology. 
Let $x\in C$ and let $\epsilon>0$ be given. Then by the previous paragraph, 
for each $P\in hull(I)$ there is a subset $V_P$ of $Prime(C)$ which is open in 
the relative Scott topology and such that $\Vert x+Q\Vert <\Vert 
x+P\Vert+\epsilon$ for all $Q\in V_P$. Set $V=\bigcup\{ V_P:P\in hull(I)\}$. 
Then $V$ is a Scott open subset of $Prime(C)$ containing $hull(I)$. Thus by 
Lemma 3.1, eventually $hull(I_{\alpha})\subseteq V$, implying that eventually 
$\Vert x+I_{\alpha}\Vert^u<\Vert x+I\Vert^u+\epsilon$. Thus $\tau_{uc}$ is 
weaker than the Scott topology on $Id^u(C)$. 

To show that $\tau_{uc}$ is stronger than the Scott topology on $Id^u(C)$, let 
$I, J\in Id^u(C)$ with $I\not\subseteq J$. Let $x\in I\backslash J$, with 
$\Vert x+J\Vert^u=1$. Then the sets $$\{ K\in Id^u(C):\Vert x+K\Vert^u 
>1/2\}$$ and 
$$\{ K\in Id^u(C):\Vert x+K\Vert^u <1/2\}$$ are disjoint neighbourhoods of $J$ 
and $I$, open in the lower and $\tau_{uc}$ topologies respectively. Thus the 
bitopological space $(Id^u(C), {\rm lower}, \tau_{uc})$ is `pseudo-Hausdorff' 
in the sense of [\the\Ko]. On the other hand, the bitopological space 
$(Id^u(C), {\rm lower}, {\rm Scott})$ is `stable', being `joincompact' (the 
definitions of these various properties are given in [\the\Ko]), so [\the\Ko; 
3.4(d)] shows that $\tau_{uc}$ is stronger than the Scott topology. Thus 
$\tau_{uc}$ and the Scott topology coincide. 

Finally, since $\tau_c$ is the join of the $\tau_{lc}$ and $\tau_{uc}$ 
topologies, while the Lawson topology is the join of the lower and Scott 
topologies, it follows immediately that $\tau_c$ coincides with the Lawson 
topology. Q.E.D. 
\bigskip 
\noindent {\bf Lemma 3.3} {\sl Let $A$ and $B$ be C$^*$-algebras with Haagerup 
tensor product $A\otimes_hB$. Let $I,J\in Id(A\otimes_h B)$ with $J^u\ll I^u$ 
in $Id^u(A\otimes_h B)$. Then $J\subseteq I_l$.} 
\bigskip 
\noindent {\bf Proof.} Let $( I_{\alpha})$ be the net of finite sums of 
product ideals contained in $I_l$, upward directed by inclusion. Then 
$(I^u_{\alpha})$ is a directed set in $Id^u(A\otimes_h B)$. If $P$ is any 
primitive ideal of $A\otimes_h B$ not containing $I$, then there is a product 
ideal contained in $I$ but not in $P$ (see the proof of Theorem 2.7). It 
follows that $I^u_{\alpha}\nearrow I^u$ in the lattice $Id^u(A\otimes_h B)$. 
Thus eventually $I^u_{\alpha}\supseteq J^u\supseteq J$, by assumption. But 
each $I_{\alpha}$ has spectral synthesis [\the\ASS; 6.4], so 
$I_{\alpha}=I^u_{\alpha}$. Thus eventually $I_{\alpha}\supseteq J$. Since 
$I_l$ is the norm closure of the product ideals which it contains [\the\ASS; 
6.7(i)], we have shown that $I_l\supseteq J$. Q.E.D. 
\bigskip 
\noindent We now need some information about quotient norms. 
\bigskip 
\noindent {\bf Lemma 3.4} {\sl Let $A$ and $B$ be C$^*$-algebras. Let $\{ 
I_i\}$ and $\{J_i\}$ $(1\le i\le n)$ be finite sets of closed ideals of $A$ 
and $B$ respectively such that all the sets $\{ P\in Prim(A): P\supseteq 
I_i\}$ $(1\le i\le n)$ are disjoint in $Prim(A)$ and all the sets $\{ Q\in 
Prim(B): Q\supseteq J_i\}$ $(1\le i\le n)$ are disjoint in $Prim(B)$. Set 
$K=\bigcap_{i=1}^n \left(I_i\otimes_h B+A\otimes_h J_i\right)$. Then for $x\in 
A\otimes_h B$, $$\Vert x+K\Vert=\max\{\Vert x+(I_i\otimes_h B+A\otimes_h 
J_i)\Vert: 1\le i\le n\}.$$} 
\noindent {\bf Proof.} Set $G=\bigcap_{i=1}^n I_i$ and $H=\bigcap_{i=1}^n 
J_i$. Then $A\otimes_hB/(G\otimes_h B+A\otimes_hH)$ is isometrically 
isomorphic to $A/G\otimes_h B/H$, see [\the\ASS; 2.6]. But 
$A/G\cong\bigoplus_{i=1}^n A/I_i$ and $B/H\cong\bigoplus_{i=1}^n B/J_i$, so 
$A/G\otimes_h B/H$ is isometrically isomorphic to $\bigoplus_{i,j=1}^n 
A/I_i\otimes_h B/J_j$. Hence $$\Vert x+(G\otimes_h B+A\otimes_h H)\Vert 
=\max\{\Vert x+(I_i\otimes_h B+A\otimes_h J_j)\Vert: 1\le i,j\le n\}.$$ Thus 
$A\otimes_h B/K$ is isometrically isomorphic to $\bigoplus_{i=1}^n 
A/I_i\otimes_h B/J_i$ and $$\Vert x+K\Vert=\max\{\Vert x+(I_i\otimes_h 
B+A\otimes_h J_i)\Vert: 1\le i\le n\}.$$ Q.E.D. 
\bigskip 
\noindent In fact we only use Lemma 3.4 in the case when the ideals $I_i$ and 
$J_i$ are maximal ideals of $A$ and $B$ respectively. We state this case 
separately. 
\bigskip 
\noindent {\bf Proposition 3.5} {\sl Let $A$ and $B$ be C$^*$-algebras and let 
$F=\{ J_1,\ldots, J_n\}$ be a finite set of maximal ideals of $A\otimes_h B$. 
Set $I(F)=\bigcap_{i=1}^n J_i$. Then for each $x\in A\otimes_h B$, $$\Vert 
x+I(F)\Vert=\max\{\Vert x+J_i\Vert: 1\le i\le n\}.$$} 
\noindent{\bf Proof.} For each $J_i\in F$, there exist maximal ideals $M_i$ of 
$A$ and $N_i$ of $B$ such that $J_i=M_i\otimes_h B+A\otimes_h N_i$ [\the\ASS; 
5.6]. The sets $\{ P\in Prim(A): P\supseteq M_i\}$ $(1\le i\le n)$ and $\{ 
Q\in Prim(B): Q\supseteq N_i\}$ $(1\le i\le n)$ are trivially disjoint in 
$Prim(A)$ and $Prim(B)$ respectively. Hence the result follows from Lemma 3.4. 
Q.E.D. 
\bigskip 
\noindent In the case when $A$ and $B$ are commutative C$^*$-algebras, 
Proposition 3.5 shows that every finite subset of $Max(A\otimes_h B)$ is a 
Helson set with Helson constant $1$. Thus $A\otimes_h B$ has the Helson 
property with constant $1$, so it follows from Theorem 1.9 that 
$A\otimes_h B$ has spectral synthesis if and only if $\tau_r$ is Hausdorff. 
Our final result, the main one of this section, is a partial extension of this 
to the non-commutative situation. The proof is similar to that of Theorem 1.9 
but with a few extra turns. 
\bigskip 
\noindent {\bf Theorem 3.6} {\sl Let $A$ and $B$ be C$^*$-algebras with 
$Prime(A)$ and $Prime(B)$ compact $T_1$-spaces. Then the following are 
equivalent: 

(i) $A\otimes_h B$ has spectral synthesis, 

(ii) $Id(A\otimes_h B)$ is $\tau_r$-Hausdorff. 
} 
\bigskip 
\noindent {\bf Proof.} (i)$\Rightarrow$(ii) follows from Theorem 2.7 and 
[\the\Id; 3.1.1]. 

(ii)$\Rightarrow$(i). Suppose that spectral synthesis fails. We shall show 
that $\tau_r$ is not Hausdorff on $Id(A\otimes_h B)$. The hypotheses on 
$Prime(A)$ and $Prime(B)$ imply that $Prime(A\otimes_h B)=Max(A\otimes_h B)$. 
Let $I\in Id(A\otimes_h B)$ with $I_l\ne I^u$. Let $N=hull(I^u)$, and let 
$(N'_{\alpha})_{\alpha}$ be the net of decreasing, Scott open subsets of 
$Max(A\otimes_h B)$ containing $N$. Let 
$(N_{\alpha})_{\alpha}$ be the net of decreasing, hull-kernel closed subsets 
of $Max(A\otimes_h B)$, where for each $\alpha$, $N_{\alpha}$ is the 
hull-kernel closure of $N'_{\alpha}$. Since $Max(A\otimes_h B)$ is locally 
compact in the hull-kernel topology, and a subset of $Max(A\otimes_h B)$ is 
open in the Scott topology if and only its complement is compact and 
saturated, it follows that $\bigcap_{\alpha} N_{\alpha}=N$. 

For a hull-kernel closed subset $M$ of $Max(A\otimes_h B)$, let 
$I(M)=\bigcap\{ Q: Q\in M\}$. Then $(I(N_{\alpha}))_{\alpha}$ is an increasing 
net in $Id^u(A\otimes_h B)$, and each $I(N_{\alpha})\ll I^u$ by the remark 
before Lemma 3.1. Hence $I(N_{\alpha})\subseteq I_l$ for each $\alpha$, by 
Lemma 3.3. Thus $J:=\overline{\bigcup_{\alpha} I(N_{\alpha})}\subseteq I_l$. 
On the other hand, $I_l$ is the smallest closed ideal with hull equal to $N$ 
[\the\ASS; 6.6(iv)], so $J=I_l$. 
For each $\alpha$, let $(F_{\beta(\alpha)})_{\beta(\alpha)}$ be the 
increasing net of finite subsets of $N_{\alpha}$. Then 
$(I(F_{\beta(\alpha)}))_{\beta(\alpha)}$ is a decreasing net in 
$Id(A\otimes_h B)$, and $\bigcap_{\beta(\alpha)}I({F_{\beta(\alpha)}}) = 
I(N_{\alpha})$. Hence $$I({F_{\beta(\alpha)}})\mapright{\beta(\alpha)} 
I(N_{\alpha})\ \ \ \ \ \ (\tau_r)$$ in $Id(A\otimes_h B)$ by Lemma 0.1. But 
$I(N_{\alpha})\to 
I_l$ $(\tau_r)$, also by Lemma 0.1, so if $(I(F_{\gamma}))_{\gamma}$ denotes 
the 
`diagonal' net, see [\the\Kel; \S 2, Theorem 4], then $I(F_{\gamma})\to I_l$ 
$(\tau_r)$ in $Id(A\otimes_h B)$. 
Now suppose that $x\notin I^u$. Then there exists $P\in hull(I)$ such that 
$x\notin P$. For each $\alpha$, eventually $P\in F_{\beta(\alpha)}$, so 
eventually $$\Vert x+I(F_{\beta(\alpha)})\Vert\ge \Vert x+P\Vert>0.$$ Hence 
the `diagonal' net $I(F_{\gamma})\to I^u$ $(\tau_n)$. 
On the other hand, if $x\in I^u$ and $\epsilon>0$ then by Proposition 3.2 
there exists a Scott open subset $M$ of $Prim(A\otimes_h B)=Max(A\otimes_h B)$ 
containing $N$ such that $\Vert x+Q\Vert<\epsilon$ for all $Q\in M$. Since 
$Max(A\otimes_h B)\backslash M$ is compact and $\bigcap N_{\alpha}=N$, there 
exists $\alpha_0$ such that for $\alpha\ge\alpha_0$, 
$N_{\alpha}\subseteq M$. Hence for $\alpha\ge\alpha_0$, $$\Vert 
x+I(F_{\beta(\alpha)})\Vert<\epsilon$$ for all $\beta(\alpha)$, by Proposition 
3.5. Hence $I(F_{\gamma})\to I^u$ $(\tau_u)$, 
using [\the\Id; 2.1], so 
$$I(F_{\gamma})\to I^u~~(\tau_r).$$ 
Since $I_l\ne I^u$, $\tau_r$ is not Hausdorff. Q.E.D. 
\bigskip 
\noindent Theorem 3.6 is sufficiently general to make one think that the 
result probably holds true for all C$^*$-algebras $A$ and $B$. 
\bigskip 
\bigskip 
\centerline{\bf References} 
\medskip 
\input TempReferences.tex 
\bigskip 
\bigskip 
\centerline{School of Mathematical Sciences} 
\centerline{University of Nottingham} 
\centerline{ NG7 2RD} 
\centerline{ U.K.} 
\medskip \centerline{email: Joel.Feinstein@nottingham.ac.uk} 
\bigskip 
\centerline{Department of Mathematical Sciences} 
\centerline{University of Aberdeen} 
\centerline{AB24 3UE} 
\centerline{U.K.} 
\medskip 
\centerline{e-mail: ds@maths.abdn.ac.uk} 
\end